\documentclass[twocolumn]{autart-1}    % Enable this line and disable the
\usepackage{xcolor}
\definecolor{darked}{RGB}{139,0,0}

\newcommand{\eproof}{\hfill\rule{2mm}{2mm}}
\newcommand{\bstate}{\begin{state} }
\newcommand{\estate}{ \hfill  \rule{1mm}{2mm} \end{state}}

\newcommand{\bass}{\begin{ass} }
\newcommand{\eass}{ \hfill  \rule{1mm}{2mm} \end{ass}}

\newcommand{\bpro}{\begin{property}}
\newcommand{\epro}{\hfill\rule{1mm}{2mm}\end{property}}

\newcommand{\brem}{ \begin{remark}  }
\newcommand{\erem}{\hfill \rule{1mm}{2mm}
\end{remark} }
\newcommand{\bthm}{\begin{theorem}  }
\newcommand{\ethm}{ \hfill  \rule{1mm}{2mm}
\end{theorem} }
\newcommand{\blem}{\begin{lemma}  }
\newcommand{\elem}{ \hfill \rule{1mm}{2mm}
\end{lemma} }
\newcommand{\bcorollary}{\begin{corollary}  }
\newcommand{\ecorollary}{  \hfill \rule{1mm}{2mm}
\end{corollary} }
\newcommand{\bdefn}{\begin{definition}}
\newcommand{\edefn}{  \hfill \rule{1mm}{2mm}
\end{definition} }
\newcommand{\bproposition}{\begin{proposition} }
\newcommand{\eproposition}{\hfill \rule{1mm}{2mm}
\end{proposition} }
\newcommand{\bexample}{\begin{example} \rm}
	\newcommand{\eexample}{ \hfill \rule{1mm}{2mm}
\end{example} }

\newcommand{\bcon}{\begin{condition} \rm}
\newcommand{\econ}{ \hfill \rule{1mm}{2mm}
\end{condition} }

\newcommand{\proofnow}{{\bf Proof: }}

\newtheorem{theorem}{\bf Theorem}[section]
\newtheorem{ass}{\bf Assumption}[section]
\newtheorem{lemma}{\bf Lemma}[section]
\newtheorem{definition}{\bf Definition}[section]
\newtheorem{remark}{\bf Remark}[section]
\newtheorem{corollary}{\bf Corollary}[section]
\newtheorem{proposition}{\bf Proposition}[section]
\newtheorem{example}{\bf Example}[section]
\newtheorem{condition}{\bf Condition}[section]
\newtheorem{state}{\bf Assumption}[section]
\newtheorem{property}{\bf Property}[section]

\usepackage{amsmath,amssymb}
\usepackage{graphicx}
\usepackage{epstopdf}
\usepackage{float}
\usepackage{multirow}%multirow table
\usepackage{cite}
\usepackage{subfigure}
\usepackage{xcolor}
\usepackage{amsfonts}
\usepackage{mathrsfs}
\usepackage{pifont}
\usepackage{stfloats}
\usepackage{color}
\usepackage{natbib}
\usepackage{enumitem}
\newcommand{\sint}{\textstyle{\int}}
\newcommand{\ssum}{\textstyle{\sum}}
\newcommand{\diag}{\mbox{diag}}
\renewcommand{\t}{^{\mbox{\tiny {T}}}}

\allowdisplaybreaks

\usepackage{graphicx}          % Include this line if your
\renewcommand{\t}{^{\ensuremath{\scriptscriptstyle \top}}}

\begin{document}

\begin{frontmatter}

\title{Small-Gain Theorem Based Distributed Prescribed-Time Convex Optimization
For Networked Euler-Lagrange Systems\thanksref{footnoteinfo}} % Title, preferably not more
\vspace{-2em}
\thanks[footnoteinfo]{This paper was not presented at any IFAC
meeting. %Corresponding author Lijun Zhu.
}
\author[China1]{Gewei Zuo}\ead{gwzuo@hust.edu.cn},
\author[Japan]{Mengmou Li}\ead{mmli,research@gmail.com},    % Add the
\author[China1,China2]{Lijun Zhu}\ead{lijun@hust.edu.cn},               % e-mail address
 % (ead) as shown

\address[China1]{School of Artificial Intelligence and Automation, Huazhong University of Science and Technology, Wuhan 430074, China}  % Please supply
\address[China2]{Key Laboratory of Imaging Processing and Intelligence Control, Huazhong University of Science and Technology, Wuhan 430074, China}             % full addresses
\address[Japan]{Graduate school of advanced science and engineering, Hiroshima University, Higashi-Hiroshima 739-0046, Japan}        % here.

                             % keyword list or with the
                                          % help of the Automatica
                                          % keyword wizard

\begin{abstract}                          % Abstract of not more than 200 words.
In this paper, we address the distributed prescribed-time convex optimization
(DPTCO) for a class of networked Euler-Lagrange systems under undirected
connected graphs. By utilizing position-dependent measured gradient
value of local objective function and local information interactions among
neighboring agents, a set of auxiliary systems is constructed to
cooperatively seek the optimal solution. The DPTCO problem is then
converted to the prescribed-time stabilization problem of
an interconnected
error system. A prescribed-time small-gain criterion is proposed to
characterize prescribed-time stabilization of the
system,
offering a novel approach that enhances the effectiveness beyond existing
asymptotic or finite-time stabilization of
an interconnected system. Under the
criterion and auxiliary systems, innovative  adaptive prescribed-time local
tracking controllers are designed for subsystems. The prescribed-time
convergence lies in the introduction of time-varying gains which increase
to infinity as time tends to the prescribed time. Lyapunov function
together with prescribed-time mapping are used to prove the prescribed-time
stability of closed-loop system as well as the boundedness of internal
signals. Finally, theoretical results are verified by one numerical
example.
\end{abstract}

\end{frontmatter}

\section{Introduction}
Cooperative control for networked Euler-Lagrange systems (NELSs) has
attracted significant attention and made a lot of progress, for example,
consensus in \cite{nuno2013consensus12,abdessameud2018consensus13},
formation control in \cite{viel2022distributed14,viel2019distributed15},
leader-follower control in \cite{abdessameud2016leader16,lu2019leader17}
and containment control in \cite{li2018two18}. Recently, another fundamentally important issue called distributed convex optimization (DCO)
arose
in cooperative control for NELSs, which
investigates
how the robots in a network can cooperatively work
to solve an optimization
problem
% cooperatively
in a distributed manner.
% where
In a DCO problem, the global objective function is composed of a sum of local objective
% function,
functions,
each of which
is known to only one agent. The existing literature on DCO can be found in  \cite{kia2015distributed03,lin2016distributed20,nedic2010constrained32,rahili2016distributed33,Gong34} and references therein.
Additionally, as a combination of distributed cooperative control and optimization, DCO has been applied to NELSs to enable multi-robots to complete practical tasks.
In \cite{zhang2017distributed05}, the feedback linearization method is used to solve the DCO problem for a class of heterogeneous NELSs without uncertainties. It is proved that the position of each subsystem achieves exponential convergence towards the optimal solution of the global objective function. To handle uncertainties, \cite{zhang2017distributed05} further proposes an adaptive optimization algorithm in which optimal trajectory generators are introduced into the feedback loop to achieve asymptotic convergence.
The adaptive optimization algorithm for NELSs
has also been considered in \cite{zou2020adaptive09} which removes
the requirements on the exact information of Lipschitz constant and
strongly convex constant of local objective functions. In \cite{zou2021distributed10},
distributed constrained convex optimization for NELSs is considered,
where the position of each robot achieves an agreement to the optimal solution and remains in the predefined closed convex constraint via projection algorithm and zero-gradient-sum algorithm.
In \cite{an2021collisions11},
by co-design optimal coordination strategy and collision avoidance
mechanism, a secure DCO algorithm is proposed for NELSs to achieve
optimization objective while avoiding collisions with other robots.

Since the relative degree of the Euler-Lagrange system (ELS) is greater
than one, the results in \cite{zhang2017distributed05,zou2020adaptive09,zou2021distributed10,an2021collisions11}
contain
an optimal trajectory generator design, which is, constructing
auxiliary systems to cooperatively seek the optimal solution and the
seeking trajectories are considered as reference trajectories for
local agents. Then the DCO problem is divided into two parts: distributed
optimal solution seeking and local reference trajectory tracking.
Indeed, the construction of auxiliary systems in \cite{zhang2017distributed05,zou2020adaptive09,zou2021distributed10,an2021collisions11}
necessitates the availability of gradient functions of local objective functions.
This approach is feedforward (open-loop) optimization and thus the
DCO implementation rests on the stabilization of a cascaded system.
Compared with feedforward optimization, feedback-based optimization
in \cite{liu2021distributed06,qin2023distributed07} is more practicable
since it utilizes measured gradient value rather than assuming the gradient
function is known. In \cite{qin2022distributed08}, feedback-based
DCO for NELSs is considered, where
position-dependent measured gradient value is utilized to construct the auxiliary systems.
The resulting closed-loop system
features an inner-outer-loop structure, achieving
exponential stability.

In this paper, we address the distributed prescribed-time convex optimization (DPTCO) problem for a class of NELSs with parameter uncertainties. Prescribed-time control is proposed to solve the problem existing in traditional finite-time and fixed-time control, where the settling time depends on the initial state and design parameters ~\citep{song2017time19}. The existing literature on finite-time or fixed-time DCO includes \cite{lin2016distributed20,wang2020distributed21,wang2020distributed22, liu2021continuous23,chen2018fixed24,chen2021fixed25,Gong34}
and references therein.  Compared with the existing literature on DCO for NELSs ~\citep{zhang2017distributed05,zou2020adaptive09,zou2021distributed10,qin2022distributed08}, our proposed DPTCO algorithm ensures that optimization errors achieve prescribed-time convergence to zero rather than merely exponential convergence. However, this introduces the challenge that the time-varying gain inherent in prescribed-time control may destabilize the closed-loop system, especially in the presence of uncertainties.
The main contributions of this paper are summarized as follows.

First, auxiliary systems are constructed to cooperatively seek the optimal solution of global objective function by utilizing position-dependent measured gradient value and local information interactions. Compared with
% the
results in \cite{liu2021distributed06,qin2023distributed07,qin2022distributed08},
the only information transmitted among neighboring agents is the position,
which
effectively reduces communication burden and
improves robustness of the network.
However, since the gradient is only available when the position is known, this real-time measurement of gradient information results in a strong coupling between the auxiliary systems and the controlled NELSs, leading to a challenging problem.
To address this,
a set of coordination transformations are introduced and then the DPTCO problem is converted
into the prescribed-time stabilization problem of
an
interconnected error system.

Second, we propose a novel prescribed-time small-gain criterion to characterize the stabilization of an interconnected system within a prescribed time. This  criterion not only encompasses classical small-gain stability conditions for interconnected systems but also incorporates specific time-varying conditions in the sense of prescribed-time stability. It offers significant benefits over traditional small-gain theorems used for asymptotic and finite-time stabilization.
Under the criterion and auxiliary systems, the adaptive prescribed-time local tracking controller is constructed for each subsystem where the adaptive estimation is introduced to compensate for the parameter uncertainties.
Consequently, the position of each subsystem
achieves prescribed-time convergence towards the optimal solution and
maintain this optimality thereafter.
Unlike finite-time or fixed-time DCO, our approach allows the settling time to be specified \textit{a priori}, independent of initial states and control parameters, ensuring uniformity across all subsystems.

Third, as a common problem in prescribed-time control, it is difficult
to analyze the boundedness of internal signals in the closed-loop system
when there
is an adaptive estimation. We introduce the prescribed-time mapping to facilitate stability analysis and further exploit prescribed-time
convergence rate. As a result, it is proved that all internal signals
in the closed-loop systems are uniformly bounded over the global time
interval.

The rest of the paper is organized as follows. 
Section \ref{sec:Notations-and-Problem} gives the notations and preliminaries. 
Section \ref{sec:Coordinate-transformation-and} introduces the DPTCO Algorithm Design and Coordinate Transformations. 
Section \ref{sec:Internected} shows the Interconnected Error System and Prescribed-Time Small-Gain criterion.
Section \ref{sec:DPTCO-Implementation} elaborates the DPTCO implementation. 
The numerical simulation is conducted in
Section \ref{sec:Simulation} and the paper is concluded in Section \ref{sec:Conclusion}.

\section{Notations and Preliminaries\label{sec:Notations-and-Problem}}

\subsection{Notations}

$\mathbb{R}$, $\mathbb{R}_{+}$ and $\mathbb{R}^{n}$ denote the
set of real numbers, the set of non-negative real numbers, and the
$n$-dimensional Euclidean space, respectively. $t_{0}$ denotes the
initial time, $T$ the prescribed-time scale, and $\mathcal{T}_{p}:=\{t:t_{0}\leq t<T+t_{0}\}$,
$\tilde{\mathcal{T}}_{p}:=\left\{ t:T+t_{0}\leq t<\infty\right\} $
the corresponding time intervals. For $\mathcal{K}_{\infty}$ functions
$\alpha(s)$ and $\alpha'(s)$, $\alpha(s)=\mathcal{O}[\alpha'(s)]$
means that $\sup_{s\in\mathbb{R}_{+}}[\alpha'(s)/\alpha(s)]<\infty$.
The symbol $1_{N}\in\mathbb{R}^{N}$ (or $0_{N}\in\mathbb{R}^{N}$)
denotes an $N$-dimensional column vector whose elements are all $1$
(or $0$).

An undirected graph is denoted as $\mathcal{G}=(\mathcal{V},\mathcal{E})$,
where $\mathcal{V}=\{1,\cdots,N\}$ is the node set and $\mathcal{E}\subseteq\mathcal{V}\times\mathcal{V}$
is the edge set. The existence of an edge $(i,j)\in\mathcal{E}$ means
that nodes $i$, $j$ can communicate with each other. Denote by $\mathcal{A}=[a_{ij}]\in\mathbb{R}^{N\times N}$
the weighted adjacency matrix, where $(j,i)\in\mathcal{E}\Leftrightarrow a_{ij}>0$
and $a_{ij}=0$ otherwise. A self edge is not allowed, i.e., $a_{ii}=0$.
The Laplacian matrix $\mathcal{L}$ of graph $\mathcal{G}$ is denoted
as $\mathcal{L}=[l_{ij}]\in\mathbb{R}^{N\times N}$, where $l_{ii}=\ssum_{j=1}^{N}a_{ij}$,
$l_{ij}=-a_{ij}$ with $i\neq j$.
The set of neighbors of node $i$ is denoted by $\mathcal N_i = \{j\in \mathcal V\mid (i,j)\in \mathcal E\}$.
 We denote the eigenvalues of $\mathcal L$ by $\lambda_1,\cdots,\lambda_N$.
If $\mathcal{G}$ is connected, zero is a simple eigenvalue of $\mathcal L$ with the eigenvector spanned by $1_{N}$, and all the other $N-1$ eigenvalues of $\mathcal{L}$ are strictly positive. In this case, we order the spectrum of $\mathcal L$ as $\lambda_1 = 0<\lambda_2\leq  \lambda_3\leq \cdots \leq \lambda_N$.
\subsection{Problem Formulation}

We consider a class of NELSs, for $i\in\mathcal{V}$,

\begin{equation}
\begin{gathered}M^{i}(q^{i})\ddot{q}^{i}+C^{i}(q^{i},\dot{q}^{i})\dot{q}^{i}=\tau^{i},\\
y^{i}=q^{i}
\end{gathered}
\label{eq:Euler-Lag}
\end{equation}
where $q^{i},\dot{q}^{i},\ddot{q}^{i}\in\mathbb{R}^{n}$ represent
the generalized position, velocity, and acceleration vectors of $i$-th
subsystem, respectively, $y^{i}$ is the measurement output, $\tau^{i}\in\mathbb{R}^{n}$
is the control input vector, $M^{i}(q^{i})\in\mathbb{R}^{n\times n}$
is the inertia matrix and $C^{i}(q^{i},\dot{q}^{i})$ is the centripetal
and Coriolis matrix.
In this paper, we suppose that each robot is moving in the horizontal direction or that the gravitational torque has already been compensated for through feedforward control. Therefore, there is no gravity component present in the dynamics.
According to \cite{lu2019adaptive01,huang2015adaptive26},
some properties for system (\ref{eq:Euler-Lag}) are presented as
follows.

\bpro \label{prop:1-1} $\dot{M}^{i}(q^{i})-2C^{i}(q^{i},\dot{q}^{i})$
is
skew-symmetric.
\epro

\bpro\label{prop:1-2} For any $x_{1},x_{2}\in\mathbb{R}^{n}$ and $i\in \mathcal V$,  $M^{i}(q^{i})x_{1}+C^{i}(q^{i},\dot{q}^{i})x_{2}=\Omega^{i}(q^{i},\dot{q}^{i},x_{1},x_{2})\theta^{i}$,
where $\Omega^{i}(q^{i},\dot{q}^{i},x_{1},x_{2})\in\mathbb{R}^{n\times p}$
is the regression matrix satisfying $\Omega^{i}(q^{i},0,0,0)=0$ and
$\theta^{i}\in\mathbb{R}^{p}$ is an unknown constant vector consists
of all inertial parameters of each subsystem in (\ref{eq:Euler-Lag}).
Furthermore, $\|\Omega^{i}(q^{i},\dot{q}^{i},x_{1},x_{2})\|\leq\rho^{i}(\|w^{i}\|_{2}+\|w^{i}\|_{2}^{2})$,
where $\rho^{i}$ is a known positive constant and $w^{i}=\left[(\dot{q}^{i})\t,x_{1}\t,x_{2}\t\right]\t$.
\epro

\bpro \label{prop:1-3} For $i\in \mathcal V$, there exists positive constants $k_{\underline{m}}^{i}$,
 $k_{\overline{m}}^{i}$
such that $k_{\underline{m}}^{i}I_{n}\leq M^{i}(q^{i})\leq
 k_{\overline{m}}^{i}I_{n}$.
\epro

In this paper, we consider the following convex optimization problem

\textcolor{black}{
\begin{equation}
 \min_{y}\;\ssum_{i=1}^{N}f^{i}(y^{i}),\quad
 \mbox{subject to}\;y^{i}=y^{j},\; \forall i\neq j
\label{eq:opti_problem}
\end{equation}
 where $y=[(y^{1})\t,\cdots,(y^{N})\t]\t$ is the lumped output, and
$f^{i}(y^{i}),i\in\mathcal{V}$ is local objective function, which is convex
and known to only agent $i$. Different from the results in \cite{zhang2017distributed05,huang2019distributed27,tang2020optimal28}
that gradient function $\nabla f^{i}(\cdot)$ of local objective function
is available, in this paper, only output-dependent measured gradient
value $\frac{\mathrm{d}f^{i}(y^{i})}{\mathrm{d}y^{i}}$ is available
for controller design and we denote it as $\nabla f^{i}(y^{i})$ for
simplicity. The global objective function $\ssum_{i=1}^{N}f^{i}(y^{i})$
is also assumed to be convex. Due to the equality constraints, the
optimal solution $y^{*}$ has the form $y^{*}=1_{N}\otimes z^{*}$
 for some $z^{*}\in\mathbb{R}^{n}$. }

\begin{figure}
\centering
\includegraphics[width=1\linewidth]{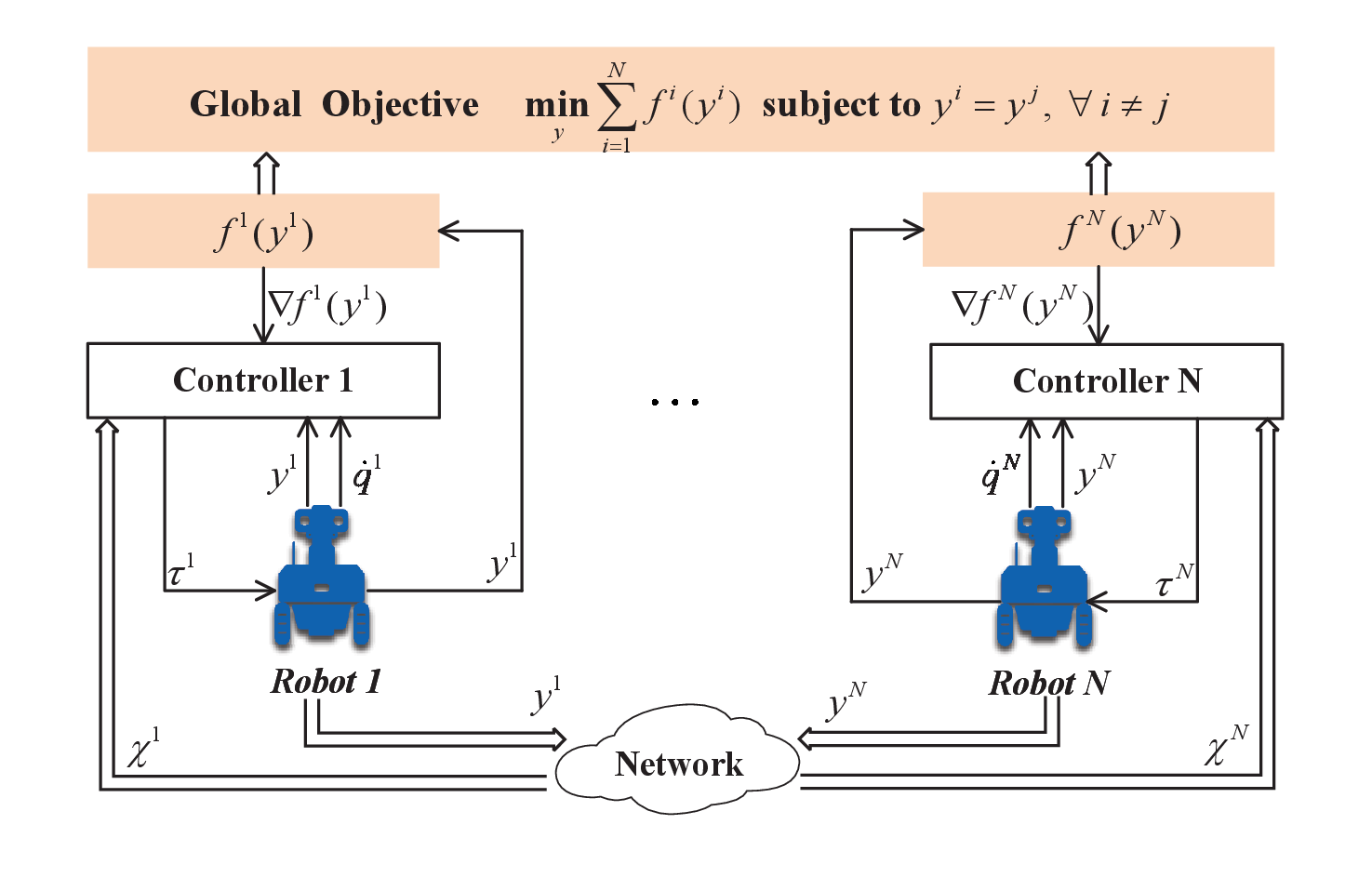}
\caption{ Schematic diagram of $N$
controlled EL systems, where $\chi^{i}=\ssum_{j\in\mathcal{N}_{i}}(y^{j}-y^{i})$
is the relative information received by $i$th agent from its neighbors.
Only system output $y^{i}$ is available for neighboring agents via
communication network.}
\label{fig:SysStructure}
\end{figure}

To make the optimization problem solvable and stability
analysis implementable, we have the following assumptions.

\bass\label{ass:graph} The undirected graph $\mathcal{G}$ is connected.\eass

\bass \label{ass:solvable} The optimal value of problem (\ref{eq:opti_problem}),
namely $f^{*}$, is finite and the optimal solution set
\[
\mathcal{\textup{Y}}_{\mbox{\small{opt} }}=\left\{ y=1_{N}\otimes z^{*}\mid\ssum_{i=1}^{N}f^{i}(z^{*})=f^{*}\right\}
\]
 is nonempty and compact ~\citep{wang2010control02}. \eass

\bass \label{ass:cost_func} For $i\in\mathcal{V}$, the objective function
$f^{i}$ is first-order differentiable, and $f^{i}$ as well as its
measured gradient $\nabla f^{i}(y^{i})$ are only known to $i$th
agent. Moreover, it is $\rho_{c}$-strongly convex and has $\varrho_{c}$-Lipschitz
gradients, i.e., for $x_{1},x_{2}\in\mathbb{R}^{n}$,
\begin{gather}(\nabla f^{i}(x_{1})-\nabla f^{i}(x_{2}))\t(x_{1}-x_{2})\geq\rho_{c} \|x_{1}-x_{2}\|^{2},\label{eq:strong-constant}\\
\|\nabla f^{i}(x_{1})-\nabla f^{i}(x_{2})\|\leq\varrho_{c} \|x_{1}-x_{2}\|\label{eq:Lipschitz-constant}
\end{gather}
where $\rho_{c}$ and $\varrho_{c}$ are some positive constants.
\eass
\brem
For $i\in \mathcal V$, if the local objective functions are defined as $f^i(y^i)= (y^i)\t p^i y^i$ where $p^i$ is a positive definite matrix, then the gradient functions are given by $\nabla f^i(y^i)= 2 p^i y^i$. In this paper, we assume that only output-dependent measured gradient values are available, specifically $2 p^i y^i$.  However, if each subsystem has \emph {a priori} knowledge of the gradient function $\nabla f^i(\cdot)$, then $\nabla f^i(z^i) =2 p^i z^i$ can be computed for any $z^i\in \mathbb R^n$. With this information, we can design $z^i$-dynamics for each $i\in \mathcal V$  to cooperatively search for the optimal solution. Since the $z^i$-dynamics are independent of the subsystems, designing local tracking controllers becomes straightforward.
 \erem

 \brem
Compared to having \emph{a priori} knowledge of the gradient function $\nabla f^i(\cdot)$, using output-dependent measured gradient values is more applicable in practical situations. For instance, in the source seeking problem ~\citep{kim2014cooperative35}, each robot cannot directly measure the distance to the source. Instead, they rely on radiation intensity sensors to estimate their proximity to the source. Similarly, in the problem of optimal formation control of multiple robots ~\citep{wu2022distributed36}, robots measure their distances from one another in real-time to determine their next direction of movement, rather than generating an optimal trajectory beforehand and then tracking it.
 \erem 

\brem Under Assumption  \ref{ass:cost_func},
the strong convexity of local objective function $f^i(y^i)$ implies its strict convexity. Therefore,
the global objective function $\ssum_{i=1}^{N}f^{i}(y^{i})$ is strictly convex, since  each  $f^i(y^i)$ is strictly convex for $i\in\mathcal V$. By Assumption \ref{ass:solvable}, the optimal solution of optimization problem \eqref{eq:opti_problem} is unique. In this paper, our proposed algorithm is also applicable when
different local objective functions have distinct strongly convex constants $\rho_c^i$ and Lipschitz constants $\varrho_c^i$. In this scenario,
inequalities \eqref{eq:strong-constant} and \eqref{eq:Lipschitz-constant}  remain valid with $\rho_c = \min\{\rho_c^1,\cdots,\rho_c^N\}$ and  $\varrho_c = \max\{\varrho_c^1,\cdots,\varrho_c^N\}$.
\erem

The objective of the DPTCO is  to design distributed controllers $u^{i},i\in\mathcal{V}$ using local output interactions $\chi^{i}=\ssum_{j\in\mathcal{N}_{i}}(y^{j}-y^{i})$ such that the outputs of all subsystems
reach the optimal solution $y^{*}$ within a prescribed-time $T+t_{0}$ and
remain at the optimal solution afterward, i.e.,
\begin{gather}
\lim_{t\to T+t_{0}}\|y(t)-y^{*}\|=0,\label{eq:objective1}\\
y(t)=y^{*},\quad t\geq T+t_{0}\label{eq:objective2}
\end{gather}
 where $y^{*}\in\textup{Y}_{\mbox{\small{opt} }}$. Furthermore, the
position $y^{i}$, velocity $\dot{q}^{i}$ and control input $\tau^{i}$
must be uniformly bounded, i.e.,
\[
\|\left[(y^{i}(t))\t,(\dot{q}^{i}(t))\t,(\tau^{i}(t))\t\right]\|\in\mathcal{L}_{\infty}
\]
 holds for $i\in\mathcal{V}$ and $t\geq t_{0}$. \textcolor{black}{The
system structure for distributed feedback optimization is shown in
Fig. \ref{fig:SysStructure}. }

\subsection{Prescribed-time Analysis}
We introduce the following time-varying functions which will be used
to derive the prescribed-time convergence.

\bdefn \label{Def:K_T} Define $\mathbb{R}_{p}=[b,\infty)$ with
$b\geq0$. A continuous differentiable function $\mu(t):\mathcal{T}_{p}\mapsto\mathbb{R}_{p}$
is said to belong to class $\mathcal{K}_{T}$ if it is strictly increasing
to infinity and
\begin{gather}
\frac{\mathrm{d}\mu(t)}{\mathrm{d}t}\Big|_{t=t_{0}}=b',\quad\lim_{t\to T+t_{0}}\frac{\mathrm{d}\mu(t)}{\mathrm{d}t}=\infty,\label{eq:mu_derivative}\\
\frac{\mathrm{d}\mu(t)}{\mathrm{d}t}\leq\tilde{b}(\mu(t))^{2},\quad\forall t\in\mathcal T _p\label{eq:bound_mu_derivative}
\end{gather}
where $b'$ and $\tilde{b}$ are some positive constant possibly related
to $T$.  \edefn

\brem
There are numerous functions
$\mu(t)$ that belong to class $\mathcal K_T$. For instance, consider $\mu(t)=(T/(T+t_0-t))^m$ with $m\geq 1$. In this case, we have  $b=1$, $b'= m/T$ and $\tilde b = m/T$. Another example is  $\mu (t)=\exp\left( 1/(T+t_0-t)\right)$, where  $b=\exp(1/T)$, $b'=\exp\left(1/T\right)(1/T)^2$, and $\tilde b = 4/\exp(2)$. It is also noteworthy that the asymptotic stability is recovered if $T$ is set to infinity.
\erem

We simplify $\mu(t)$ and $\dot{\mu}(t)$ as $\mu$ and $\dot{\mu}$
throughout this paper if no confusion
occurs and define
\begin{gather}
\kappa^{\iota}(\alpha(\mu))=\exp\left(\iota\sint_{t_{0}}^{t}\alpha(\mu(\tau))\mathrm{d}\tau\right),\;t\geq t_{0},\label{eq:kappa}\\
\kappa^{\iota}(\alpha(\mu(\tau)))=\exp\left(\iota\sint_{t_{0}}^{\tau}\alpha(\mu(s))\mathrm{d}s\right),\;\tau\geq t_{0}\label{eq:kappa_tau}
\end{gather}
 where $\iota$ is constant and $\alpha\in\mathcal{K}_{\infty}$.
We note $\kappa^{\iota}(\alpha(\mu))$ converges to zero as $t\to T+t_{0}$
for any $\iota<0$ and $\alpha\in\mathcal{K}_{\infty}$.

\bdefn A function $\alpha:[0,\infty)\mapsto[b,\infty)$ is said to
belong to class $\mathcal{K}_{\infty}^{e}$,  if it is strictly increasing
and $\alpha(0)=b>0$. \edefn

\bdefn A continuous function $\beta:[0,\infty)\times[0,\infty)\mapsto[0,\infty)$
is said to belong class $\mathcal{KL}_{T}$ if, for each fixed $s$,
the mapping $\beta(r,s)$ belongs to class $\mathcal{K}_{\infty}$
with respect to $r$ and, for each fixed $r$, there exists a constant
$T$ such that, for $s\in[0,T)$, the mapping $\beta(r,s)$ is decreasing
with respect to $s$ and satisfies $\beta(r,s)\to0$ as $s\to T$,
$\beta(r,s)=0$ for $s\in[T,\infty)$. In addition, the function $\beta$ is said
to belong class  $\mathcal{KL}_{T}^{e}$ if, for each fixed $s$,
the mapping $\beta(r,s)$ belongs to class $\mathcal{K}_{\infty}^{e}$
with respect to $r$.
\edefn

\bdefn \textcolor{black}{\label{def:PT-Stbility}Consider the system
\begin{equation}
\dot{\chi}=g(t,\chi,d(t)),\quad\chi(t_{0})=\chi_{0}\label{eq:chi-dynamics}
\end{equation}
where $\chi\in\mathbb{R}^{n}$ is the state, $d(t):[t_{0},\infty)\mapsto\mathbb{R}^{n_{d}}$
is the external input, and $\chi_{0}$ is the initial state at initial
time $t_{0}$. For any given $T>0$, the $\chi$-dynamics is said
to be prescribed-time stable if there exits $\beta\in\mathcal{KL}_{T}^{e}$
such that, for any $\chi_{0}\in\mathbb{R}^{n}$, $d\in\mathbb{R}^{n_{d}}$,
\begin{equation*}
\|\chi(t)\|\leq\beta(\|\chi_{0}\|,t-t_{0})
\end{equation*}
holds for $t_{0}\leq t<\infty$.} \edefn

\bdefn \label{defn:PTISS}The continuously differentiable function
$V(\chi):\mathbb{R}^{n}\mapsto\mathbb{R}_{+}$ is called the prescribed-time
input-to-state stable (ISS) Lyapunov function for system in \eqref{eq:chi-dynamics}
with $d$ being the input, if $V(\chi)$ and its derivative along
the trajectory of the system satisfy, for all $\chi\in\mathbb{R}^{n}$
and $t\in\mathcal{T}_{p}$,
\begin{equation}
\begin{gathered}\underline{\alpha}(\|\chi\|)\leq V(\chi)\leq  \overline{\alpha} (\|\chi\|),\\
\dot{V}\leq-\alpha(\mu)V+\tilde{\epsilon}(\mu)\epsilon(\|d\|)
\end{gathered}
\label{eq:PTISS}
\end{equation}
 where $\underline{\alpha}$, $ \overline{\alpha}$, $\alpha$, $\tilde{\epsilon}$,
$\epsilon\in\mathcal{K}_{\infty}$. \emph{  $\alpha(\mu)$, $\tilde{\epsilon}(\mu)$
and $\epsilon(\|d\|)$ are called prescribed-time convergent gain,
prescribed-time ISS gain and (normal) ISS gain, respectively.}  The
inequalities in (\ref{eq:PTISS}) are simplified as $V(\chi)\sim\{ \underline{\alpha},\overline{\alpha}, \alpha,[\epsilon(\|d\|), \tilde{\epsilon}]\mid\dot{\chi}=g(t,\chi,d)\} $.
When system \eqref{eq:chi-dynamics} has multiple inputs, i.e., $\dot \chi =g(t,\chi,d_1(t),\cdots,d_n(t))$ where $d_i \in \mathbb R^{n_{di}}$, the second inequality of (\ref{eq:PTISS})
becomes $\dot{V}\leq-\alpha(\mu)V+\ssum\tilde{\epsilon}_{i}(\mu)\epsilon_{i}(\|d_{i}\|),$
and the inequalities are simplified as $V(\chi)\sim\{\underline{\alpha}, \overline{\alpha},\alpha,[\epsilon_{1}(\|d_{1}\|),\tilde{\epsilon}_{1}], \cdots,[\epsilon_{n}(\|d_{n}\|),\tilde{\epsilon}_{n}]\mid\dot{\chi}=g(t,\chi,d_1,\cdots,d_n)\}$.\edefn

\bdefn \label{defn:mapping}
A time-varying mapping nonlinear mapping $h(\chi,\mu):\mathbb{R}^{n}\times\mathbb{R}_{p}\mapsto\mathbb{R}^{m}$
is said to be a prescribed-time mapping if the partial derivatives
$\frac{\partial h(\chi,\mu)}{\partial\chi}$ and $\frac{\partial h(\chi,\mu)}{\partial\mu}$
exist. Additionally, if $h(\chi,\mu)$ is uniformly bounded for $t\in \mathcal T_p$, this uniform boundedness guarantees the prescribed-time convergence of $\chi$.
\edefn

\brem \label{rem:mapping}The prescribed-time mapping is a very useful
tool in deriving prescribed-time stabilization of a dynamic system.
For example, consider a disturbed system $\dot{\chi}=-k\mu\chi+d(t)$,
where $\chi\in\mathbb{R}^{n}$, $k>0$, $\mu\in\mathcal{K}_{T}$, and $d(t)\in\mathbb{R}^{n}$ satisfying
$\sup_{t\geq t_{0}}\|d(t)\|\leq\bar{d}$ .
Define the Lyapunov function candidate for $\chi$-dynamics as $V(\chi)=\frac{1}{2}\chi\t\chi$,
then its time derivative satisfies $\dot{V}(\chi)=-k\mu\|\chi\|^{2}+\chi\t d(t)\leq(-2k\mu+1)V(\chi)+\bar{d}^{2}/2$.
Invoking comparison lemma yields
\begin{align*}
V(\chi(t)) & \leq\kappa^{-2k}(\mu)\exp\left(t-t_{0}\right)V(\chi(t_{0}))\\
 & \quad+\frac{\bar{d}^{2}}{2}\kappa^{-2k}(\mu)\sint_{t_{0}}^{t}\kappa^{2k}(\mu(\tau))\exp(t-\tau)\mathrm{d}\tau
\end{align*}
where $\kappa(\mu)$ and $\kappa(\mu(\tau))$ are defined in (\ref{eq:kappa})
and (\ref{eq:kappa_tau}).
Observing the upper bound of $V(\chi)$,
it is both tedious and complicated to analyze the prescribed-time convergence
of $V(\chi)$ as well as $\chi$. By applying the prescribed-time mapping,
we define $\tilde{\chi}=h(\chi,\mu)=\mu\chi$. The $\tilde{\chi}$-dynamics is given by
 $\dot{\tilde{\chi}}=(-k\mu+\dot{\mu}/\mu)\tilde{\chi}+\mu d(t)$.
We define the Lyapunov function candidate for $\tilde{\chi}$-dynamics as $\tilde{V}(\tilde{\chi})=\frac{1}{2}\tilde{\chi}\t\tilde{\chi}$
and assume $k\geq\frac{k^{*}}{2}+\tilde{b}+\frac{1}{2}$, where $\tilde{b}$
is denoted in (\ref{eq:bound_mu_derivative}). Upon using Young's
inequality, the time derivative satisfies $\dot{\tilde{V}}(\tilde{\chi})\leq-k^{*}\mu\tilde{V}(\tilde{\chi})+\frac{1}{2}\mu\bar{d}^{2}$.
Invoking comparison lemma yields $\tilde{V}(\tilde{\chi}(t))\leq\kappa^{-k^{*}}(\mu) \left(\tilde{V}(\tilde{\chi}(t_{0}))-\frac{\bar{d}^{2}}{2k^{*}}\right)+\frac{\bar{d}^{2}}{2k^{*}}
 \leq \tilde V(\tilde \chi(t_0))+\frac{\bar{d}^{2}}{2k^{*}}$,
which further implies the prescribed-time convergence of  the original
state $\chi$, i.e., $\|\chi(t)\|\leq\left(\mu(t_0)^{2}\|\chi(t_{0})\|^{2}+\bar{d}^{2}/k^{*}\right)^{\frac{1}{2}}\mu^{-1}(t)$. 
\erem

\section{DPTCO Algorithm Design and Coordinate Transformations\label{sec:Coordinate-transformation-and}}

Given that the relative degree of each subsystem in (\ref{eq:Euler-Lag}) is greater than one, directly proving whether the output of the closed-loop system can converge to the optimal solution is a very complicated problem. Instead, we construct auxiliary systems to seek the optimal solution, with these seeking trajectories considered as reference trajectories for local subsystems. The DPTCO problem is then divided into two parts: distributed optimal solution seeking and local reference trajectory tracking.

For $i\in\mathcal{V}$ and , the  auxiliary systems are designed as follows:
\begin{align}
\dot{\varpi}^{i}(t) & =-c\mu\left(\nabla f^{i}(y^{i})+\chi_y^i+v^{i}\right),\label{eq:dot_varpi_i}\\
\dot{v}^{i}(t) & =c\mu\chi_y^i,\quad t\in\mathcal{T}_{p}\label{eq:dot_v_i}
\end{align}
where $c>0$ is a design parameter, $\mu\in\mathcal{K}_{T}$ is
denoted in Definition \ref{Def:K_T}, and $\chi_y^i = \ssum_{j\in\mathcal{N}_{i}}a_{ij}(y^{i}-y^{i})$ represents the relative information received by $i$th subsystem. Additionally,
\begin{equation}
\dot{\varpi}^{i}(t)=0,\quad\dot{v}^{i}(t)=0, \quad t\in\tilde{\mathcal{T}}_{p}.\label{eq:dot_varpi_v}
\end{equation}
$\varpi^{i}\in\mathbb{R}^{n}$
aims to converge to $z^{*}$ within $T+t_{0}$ and remain as $z^{*}$
afterward, $v^{i}\in\mathbb{R}^{n}$ is designed to adaptively find
the value $\nabla f^{i}(z^{*})$.

Based on auxiliary systems (\ref{eq:dot_varpi_i})-(\ref{eq:dot_varpi_v}),
for $i\in\mathcal{V}$, the local prescribed-time tracking controller
is designed as
\begin{equation}
\tau^{i}(t)=\begin{cases}
\pi^{i}(q^{i},\dot{q}^{i},\varpi^{i},\hat{\theta}^{i},\mu), & \mbox{for}\;t\in\mathcal{T}_{p}\\
0, & \mbox{for}\;t\in\tilde{\mathcal{T}_{p}}
\end{cases}\label{eq:tau^i}
\end{equation}
where
\begin{align}
 & \pi^{i}(q^{i},\dot{q}^{i},\varpi^{i},\hat{\theta}^{i},\mu)\nonumber \\
 & =-k_{2}^{i}\mu(\dot{q}^{i}+k_{1}^{i}\mu(y^{i}-\varpi^{i})) -\Omega^{i}(q^{i},\dot{q}^{i},z_{1}^{i},z_{2}^{i})\hat{\theta}^{i}\label{eq:pi^i}
\end{align}
 and $k_{1}^{i}$, $k_{2}^{i}$ are design parameters. The matrix function $\Omega^{i}(q^{i},\dot{q}^{i},z_{1}^{i},z_{2}^{i})$
can be calculated according to Property \ref{prop:1-2},
with $z_{1}^{i}$, $z_{2}^{i}$ being intermediate variables denoted
as
\begin{align}
z_{1}^{i} & =\left((\iota-1)\tilde{\mu}+k_{1}^{i}\mu\right)\dot{q}^{i}+\iota k_{1}^{i}\mu\tilde{\mu}(y^{i}-\varpi^{i}),\label{eq:z_1^i}\\
z_{2}^{i} & =k_{1}^{i}\mu(y^{i}-\varpi^{i})\label{eq:z_2^i}
\end{align}
where $\iota>0$ is a design parameter and $\tilde{\mu}=\dot{\mu}/\mu$.
In (\ref{eq:tau^i}), $\hat{\theta}^{i}\in\mathbb{R}^{p}$ denotes
the estimate of unknown parameter $\theta^{i}$ whose dynamics is
designed as
\begin{equation}
\dot{\hat{\theta}}^{i}=\varsigma^{i}(q^{i},\dot{q}^{i},\varpi^{i},\mu)-\sigma^{i}\mu\hat{\theta}^{i} \label{eq:dot_hat_theta^i}
\end{equation}
 for $i\in\mathcal{V}$ and $t\in\mathcal{T}_{p}$, where $\sigma^{i}\geq\tilde{b}$
is  a design parameter and
\begin{eqnarray}
 & &\varsigma^{i}(q^{i},\dot{q}^{i},\varpi^{i},\mu)\nonumber \\
 & & =2\mu^{2\iota-2}\left(\Omega^{i}(q^{i},\dot{q}^{i},z_{1}^{i},z_{2}^{i})\right)\t(\dot{q}^{i}+k_{1}^{i}\mu(y^{i}-\varpi^{i})).\label{eq:varsigma^i}
\end{eqnarray}
We note that $\hat{\theta}^{i}(t)$ only exists for $t\in\mathcal{T}_{p}$.
Define
\[
\tilde{\theta}^{i}=\theta^{i}-\hat{\theta}^{i}
\]
 as the estimate error for $i\in\mathcal{V}$.

Let $y^{i}=\varpi^{i}$, $\dot{q}^{i}=0$, one has $z_{1}^{i}=0$
and $z_{2}^{i}=0$, then by Property \ref{prop:1-2}, $\pi^{i}(q^{i},\dot{q}^{i},\varpi^{i},\hat{\theta}^{i},\mu)$
satisfies $\pi^{i}(\varpi^{i},0,\varpi^{i},\hat{\theta}^{i},\mu)=0$
for any $\hat{\theta}^{i}\in\mathbb{R}^{p}$ and $\mu\in\mathbb{R}_{p}$.

Define the lumped vectors $\varpi=[(\varpi^{1})\t,\cdots,(\varpi^{N})\t]\t$,
$v=[(v^{1})\t,\cdots,(v^{N})\t]\t$, $\dot{q}=[(\dot{q}^{1})\t,\cdots,(\dot{q}^{N})\t]\t$,
$\hat{\theta}=[(\hat{\theta}^{1})\t,\cdots,(\hat{\theta}^{N})\t]\t$
and $\dot{\hat{\theta}}=[(\dot{\hat{\theta}}^{1})\t,\cdots,(\dot{\hat{\theta}}^{N})\t]\t$,
then according to (\ref{eq:Euler-Lag}), (\ref{eq:dot_varpi_i}),
(\ref{eq:dot_v_i}) and (\ref{eq:pi^i}), for $t\in\mathcal{T}_{p}$,
the closed-loop NELSs can be written in a compact form as
\begin{align}
\dot{\varpi}(t) & =-c\mu\nabla F(y)-c\mu\bar{\mathcal{L}}y-c\mu v,\label{eq:dot_varpi}\\
\dot{v}(t) & =c\mu\bar{\mathcal{L}}y,\label{eq:dot_v}\\
\dot{y}(t) & =\dot{q},\label{eq:dot_y}\\
\ddot{q}(t) & =M^{-1}(q)(\pi(q,\dot{q},\varpi,\hat{\theta},\mu)-C(q,\dot{q})\dot{q}),\label{eq:ddot_q}\\
\dot{\hat{\theta}}(t) & =\varsigma(q,\dot{q},\varpi,\mu)-\bar\sigma\mu\hat{\theta}\label{eq:dot_hat_theta}
\end{align}
 where $\bar{\mathcal{L}}=\mathcal{L}\otimes I_{n}$, $M=\mbox{diag}\{ M^{1},\cdots,M^{N}\} $,
$\pi=[(\pi^{1})\t,\cdots,(\pi^{N})]\t$ with $\pi^i$ denoted in \eqref{eq:tau^i},
$C=\mbox{diag}\{ C^{1},\cdots,C^{N}\} $, $\nabla F(y)=[(\nabla f^{1}(y^{1}))\t,\cdots,(\nabla f^{N}(y^{N}))\t]\t$, $\varsigma=[(\varsigma^{1})\t,\cdots,(\varsigma^{N})\t]\t$ and $\bar \sigma =\mbox{diag}\{\sigma^1 ,\cdots,\sigma^N\}\otimes I_p$.
Define
\begin{equation}
\eta=\left[c,\iota,k_{1}\t,k_{2}\t,\sigma\t\right]\t\in\mathbb{R}^{3N+2}\label{eq:eta}
\end{equation}
 as the vector contains all design parameters of the whole controlled
system, where $k_{1}=[k_{1}^{1},\cdots,k_{1}^{N}]\t$, $k_{2}=[k_{2}^{1},\cdots k_{2}^{N}]\t$
and $\sigma=[\sigma^{1},\cdots,\sigma^{N}]\t$.

\bproposition \label{prop-1} Let the optimal solution of optimization
problem (\ref{eq:opti_problem}) is $y^{*}=1_{N}\otimes z^{*}$ for
some $z^{*}\in\mathbb{R}^{n}$. Suppose Assumption \ref{ass:graph},
\ref{ass:solvable} hold and the initial value of $v^{i}$ satisfies
$\ssum_{i=1}^{N}v^{i}(t_{0})=0$, then the equilibrium points of system
(\ref{eq:dot_varpi})-(\ref{eq:ddot_q}) are
\begin{gather}
\varpi^{*}=1_{N}\otimes z^{*},\quad v^{*}=-\nabla F(1_{N}\otimes z^{*}),\label{eq:varpi*}\\
y^{*}=1_{N}\otimes z^{*},\quad\dot{q}^{*}=0\label{eq:y*}
\end{gather}
 where $\nabla F(1_{N}\otimes z^{*})=\left[(\nabla f^{1}(z^{*}))\t,\cdots,(\nabla f^{N}(z^{*}))\t\right]\t$.
\eproposition

\proofnow Consider the solution of
\begin{equation}
\begin{cases}0  =-c\mu\nabla F(y)-c\mu\bar{\mathcal{L}}y-c\mu v\\
0  =c\mu\bar{\mathcal{L}}y\\
0  =\dot{q}\\
0  =M^{-1}(q)(\pi(q,\dot{q},\varpi,\hat{\theta},\mu)-C(q,\dot{q})\dot{q})
\end{cases}
\label{eq:dynamic_zero}
\end{equation}
where we omit  \eqref{eq:dot_hat_theta} since the value of $\hat \theta$ at the equilibrium point does not  impact the equilibrium points of other variables.
By the second and third equations of (\ref{eq:dynamic_zero}), one
has $\dot{q}^{*}=0$ and $y^{*}=1_{N}\otimes\zeta$
 for some $\zeta\in\mathbb{R}^{n}$. Since $\mathcal{G}$ is undirected
and connected, the Laplacian matrix $\mathcal{L}$ is symmetric and
its null space is spanned by $1_{N}$, then $1_{N}\t\mathcal{L}=0$.
By $\ssum_{i=1}^{N}v^{i}(t_{0})=0$, for $t\in\mathcal{T}_{p}$,
\begin{align*}
&(1_{N}\t\otimes I_{n})v(t) \\
& =(1_{N}\t\otimes I_{n})v(t_{0})+\sint_{t_{0}}^{t}(1_{N}\t\otimes I_{n})\dot{v}(\tau)\mathrm{d}\tau\\
 & =\ssum_{i=1}^{N}v^{i}(t_{0})+\sint_{t_{0}}^{t}(1_{N}\t\mathcal{L}\otimes I_{n})c\mu(\tau)y(\tau)\mathrm{d}\tau\\
 & =0.
\end{align*}
 Left multiplying the first inequality in (\ref{eq:dynamic_zero})
by $1_{N}\t\otimes I_{n}$ yields $\ssum_{i=1}^{N}\nabla f^{i}(\zeta)=0$.
 For the optimization problem (\ref{eq:opti_problem}), the necessary
and sufficient condition for a point $1_{N}\otimes z^{*}$ to be the
optimal solution is $\ssum_{i=1}^{N}\nabla f^{i}(z^{*})=0$ ~\citep{kia2015distributed03}.
Thus we have $\zeta=z^{*}$ and $v^{*}=-\nabla F(1_{N}\otimes z^{*})$.
Note that $\pi(\varpi,0,\varpi,\hat{\theta},\mu)=0$ and $q=y$, then
the forth equality in (\ref{eq:dynamic_zero}) leads to $\varpi^{*}=1_{N}\otimes z^{*}$.
\eproof

Define error variables as
\begin{align}
e_{\varpi} & =\varpi-1_{N}\otimes z^{*},\label{eq:e_varpi}\\
e_{v} & =v+\nabla F(1_{N}\otimes z^{*}),\label{eq:v}\\
e_{y} & =y-\varpi\label{eq:e_y}
\end{align}
and
\begin{gather}
e_{r}=\left[e_{\varpi}\t,e_{v}\t\right]\t,\quad e_{s}=\left[e_{y}\t,\dot{q}\t\right]\t.\label{eq:e_r_s}
\end{gather}
Based on Proposition \ref{prop-1}, the DPTCO problem is converted
into the prescribed-time stabilization problem of $e_{r}$ and $e_{s}$.

\brem
The proof of Proposition \ref{prop-1} implies that the implementation of the DPTCO with a feedback approach requires a precondition: the controlled NELSs in \eqref{eq:Euler-Lag} must have an equilibrium point at $[(q^{i,})\t,(\dot{q}^{i,})\t]\t = [(\zeta^{i,})\t,0_n\t]\t$ for any $\zeta^{i,} \in \mathbb{R}^n$ and $i \in \mathcal{V}$. This necessitates that the equilibrium point for each robot be position-independent. Consequently, system (1) assumes that each robot is moving in the horizontal direction or that the gravitational torque has already been compensated for through feedforward control.
Furthermore, the adaptive estimation $\hat{\theta}^i$ in \eqref{eq:pi^i} is introduced to ensure that $\pi^{i}(\varpi^{i},0,\varpi^{i},\hat{\theta}^{i},\mu)=0$ for any $\hat{\theta}^{i} \in \mathbb{R}^{p}$ and $\mu \in \mathbb{R}_{p}$. This condition is crucial to extend the prescribed-time control indefinitely and ensure the continuity of the control inputs over the global time interval.
\erem

\section{Interconnected Error System and Prescribed-Time Small-Gain criterion}\label{sec:Internected}
In this section, we  show that the feedback approach based DPTCO problem, along with the boundedness issue of internals signals in the closed-loop systems,  can be converted into the practical stability problem of an interconnected system. Subsequently, we propose a prescribed-time small-gain criterion to  delineate the conditions necessary for ensuring practical stability.

\subsection{Interconnected Error System and Prescribed-Time Mapping}
To make controllers in (\ref{eq:tau^i}) implementable, the adaptive
estimate $\hat{\theta}^{i}(t)$ must be bounded for $t\in\mathcal{T}_{p}$
and $i\in\mathcal{V}$. It is
% equal prove
 equivalent to proving
the uniform boundedness of
% estimate
 the estimation
error $\tilde{\theta}^{i}$ since $\theta^{i}$ is a finite constant
vector. Furthermore, since $\mu(t)$ increases to infinity as $t\to T+t_{0}$,
the proposed design must guarantee uniform boundedness of the $\mu(t)$-dependent
terms in (\ref{eq:dot_varpi_i}), (\ref{eq:dot_v_i}), (\ref{eq:tau^i})
and (\ref{eq:dot_hat_theta^i}). This issue can be addressed by further
exploiting the prescribed-time convergence rate of $e_{r}$ and $e_{s}$.
For example, from (\ref{eq:dot_varpi}), (\ref{eq:dot_v}), and (\ref{eq:e_varpi})-(\ref{eq:e_r_s}),
one has
\begin{equation}
\|\dot{e}_{r}\|\leq c(2\lambda_{N}+\varrho_{c}+1)\mu(\|e_{r}\|+\|e_{y}\|)\label{eq:bound_dot_e_r}
\end{equation}
which implies that to guarantee the uniform boundedness of $\dot{\varpi}^{i}$
and $\dot{v}^{i}$ for $i\in\mathcal{V}$ and $t\in\mathcal{T}_{p}$,
$e_{r}$ and $e_{s}$ must achieve prescribed-time convergence rate
no less than $\mu^{-1}(t)$. As explained in Remark \ref{rem:mapping},
it is  impractical to directly analyze the prescribed-time convergence
as well as prescribed-time convergence rate of $e_{r}$ and $e_{s}$
through their dynamics. Therefore, we introduce the following time-varying state transformations, for $t\in\mathcal{T}_{p}$,
\begin{equation}
\tilde{e}_{r}=\mu^{\iota}e_{r},\quad\tilde{e}_{s}=\Lambda_{1}(\mu)e_{s}\label{eq:tilde_e_r}
\end{equation}
 where $\iota$ is introduced in (\ref{eq:z_1^i}) and
\[
\Lambda_{1}(\mu)=\left[\begin{array}{cc}
\mu^{\iota}I_{nN} & 0_{nN\times nN}\\
\mu^{\iota}\bar{K}_{1} & \mu^{\iota-1}I_{nN}
\end{array}\right]
\]
 where $\bar{K}_{1}=K_{1}\otimes I_{n}\in\mathbb{R}^{nN\times nN}$
with $K_{1}=\diag\{k_{1}^{1},\cdots,k_{1}^{N}\}$. According to (\ref{eq:Euler-Lag}), \eqref{eq:dot_hat_theta^i},
(\ref{eq:dot_varpi}) and (\ref{eq:dot_v}), for $t\in\mathcal{T}_{p}$,
the dynamics of $\tilde{e}_{r}$, $\tilde{e}_{s}$ and $\tilde{\theta}$
can be expressed as
\begin{gather}
\dot{\tilde{e}}_{r}=\zeta_{r}(\tilde{e}_{r},\tilde{e}_{s},\mu),\label{eq:dot_til_e_r}\\
\dot{\tilde{e}}_{s}=\zeta_{s}(\tilde{e}_{r},\tilde{e}_{s},\tilde{\theta},\mu),\label{eq:dot_til_e_s}\\
\dot{\tilde{\theta}}=\zeta_{\theta}(\tilde{e}_{s},\tilde{\theta},\mu)\label{eq:dot_til_theta}
\end{gather}
where
\begin{align*}
 & \zeta_{r}(\tilde{e}_{r},\tilde{e}_{s},\mu)=\iota\tilde{\mu}\tilde{e}_{r}+c\mu^{\iota+1}\\
 & \times \left[\begin{gathered}-\bar{\mathcal{L}}(e_{\varpi}+e_{y})-e_{v} - \nabla\tilde{F}(y,\varpi)- \nabla\tilde{F}(\varpi,\varpi^*)\\
\bar{\mathcal{L}}e_{y}+\bar{\mathcal{L}}e_{\varpi}
\end{gathered}
\right],\\
 & \zeta_{s}(\tilde{e}_{r},\tilde{e}_{s},\mu)=\tilde{\mu}\left[\begin{array}{cc}
\iota I_{nN} & 0_{nN\times nN}\\\bar K_{1} & (\iota-1)I_{nN}
\end{array}\right]\tilde{e}_{s}\\
 & +\Lambda_{1}(\mu)\left[\begin{array}{c}
\dot{q}-\dot{\varpi}\\
M^{-1}(q)(\pi(q,\dot{q},\varpi,\hat{\theta},\mu)-C(q,\dot{q})\dot{q})
\end{array}\right],\\
 & \begin{aligned}\zeta_{\theta}(\tilde{e}_{s},\tilde{\theta},\mu) & =\bar\sigma\mu(\theta-\tilde{\theta})\\
 & \quad-2\mu^{\iota-1}\left(\Omega(q,\dot{q},z_{1},z_{2})\right)\t\left[0,I_{nN}\right]\tilde{e}_{s}
\end{aligned}
\end{align*}
 where $\Omega=\mbox{diag}\left\{\Omega^1,\cdots,\Omega^N\right\}$ and
 we used the facts $\bar{\mathcal{L}}(1_{N}\otimes z^{*})=0$
and
\begin{align}
\dot{e}_{\varpi}  &=-c\mu\nabla F(y)-c\alpha(\mu)\bar{\mathcal{L}}y-c\alpha(\mu)v\nonumber \\
  &=-c\mu\big(\bar{\mathcal{L}}e_{\varpi}+\bar{\mathcal{L}}e_{y}+e_{v}+\nabla\tilde{F}(y,\varpi) \notag\\
  &\quad + \nabla\tilde{F}(\varpi,\varpi^*)\big )\label{eq:dot_e_varpi}
\end{align}
where $ \nabla\tilde{F}(y,\varpi)$ and $ \nabla\tilde{F}(\varpi,\varpi^*)$ are
denoted by
\begin{align*}
 \nabla\tilde{F}(y,\varpi)= & \left[(\nabla f^{1}(y^{1})-\nabla f^{1}(\varpi^{1}))\t,\right.\\
 & \left.\cdots,(\nabla f^{N}(y^{N})-\nabla f^{N}(\varpi^{N}))\t\right]\t,\\
 \nabla\tilde{F}(\varpi,\varpi^* )= & \left[(\nabla f^{1}(\varpi^{1})-\nabla f^{1}(z^{*}))\t,\right.\\
 & \left.\cdots,(\nabla f^{N}(\varpi^{N})-\nabla f^{N}(z^{*}))\t\right]\t.
\end{align*}

By examining \eqref{eq:dot_til_e_r}-\eqref{eq:dot_til_theta}, we observe that the dynamics of $\tilde{e}_{r}$, $\tilde{e}_{s}$ and $\tilde \theta$ are in  an interconnected form.

We have the following lemma, whose proof is given in appendix.
\blem \label{lem:PT-Mapping}
The time-varying state transformation
(\ref{eq:tilde_e_r}) serves as  a prescribed-time mapping from $e_{r}$,
$e_{s}$ to $\tilde{e}_{r}$, $\tilde{e}_{s}$ if $\iota>2$. For
$t\in\mathcal{T}_{p}$, the uniform boundedness of $\tilde{e}_{r}$ and $\tilde{e}_{s}$
imply the prescribed-time convergence of $e_{r}$ and $e_{s}$ with
the upper bounds
\begin{equation}
\|e_{r}\|\leq\mu^{-\iota}\gamma_{r}(\|\tilde{e}_{r}\|_{\mathcal{T}}), \quad\|e_{s}\|\leq\mu^{-\iota+1}\gamma_{s}(\|\tilde{e}_{s}\|_{\mathcal{T}})\label{eq:bounds_e_r_s}
\end{equation}
 where $\|\tilde{e}_{r}\|_{\mathcal{T}}=\sup_{t\in\mathcal{T}_{p}}\|\tilde{e}_{r}\|$,
$\|\tilde{e}_{s}\|_{\mathcal{T}}=\sup_{t\in\mathcal{T}_{p}}\|\tilde{e}_{s}\|$
and $\gamma_{r},\gamma_{s}\in\mathcal{K}_{\infty}$. Furthermore,
if $\tilde{e}_{r}$ and $\tilde{e}_{s}$ are uniformly bounded, then
$\dot{\varpi}^{i}$, $\dot{v}^{i}$ in (\ref{eq:dot_varpi}) and (\ref{eq:dot_v}),
controller $\tau^{i}$ in (\ref{eq:tau^i}), and adaptive estimate
$\hat{\theta}^{i}$ as well as its dynamics $\dot{\hat{\theta}}^{i}$
in (\ref{eq:dot_hat_theta^i}) are uniformly bounded for $t\in\mathcal{T}_{p}$
and $i\in\mathcal{V}$. \elem

\subsection{Prescribed-time small-gain criterion }
Based on Lemma \ref{lem:PT-Mapping}, the DPTCO problem is now converted into the uniform boundedness problem of $\tilde e_r $ and $\tilde e_s$.
To facilitate the prescribed-time stability analysis and provide an analytical method for selecting control parameters, we introduce the following prescribed-time small-gain criterion, whose proof is given in appendix.
\bthm \label{lem:SG-Criterion} Consider the system
\begin{equation}
\begin{gathered}\dot{\chi}_{1}=h_{1}(t,\chi_{1},\chi_{2},d_{1}(t)),\\
\dot{\chi}_{2}=h_{2}(t,\chi_{1},\chi_{2},d_{2}(t))
\end{gathered}
\label{eq:dot_chi}
\end{equation}
 where $\chi=[\chi_{1}\t,\chi_{2}\t]\t\in\mathbb{R}^{n_{1}}\times\mathbb{R}^{n_{2}}$
is the state with $\chi_{0}=[\chi_{1}\t(t_{0}),\chi_{2}\t(t_{0})]\t$
as the initial value, and $d_{i}(t):[t_{0},\infty)\mapsto\mathbb{D}_{i}$
is the external disturbance for $i=1,2$, where $\mathbb{D}_{i}$
is a compact set belonging to $\mathbb{R}^{l_{i}}$. Each $\chi_{i}$-subsystem
admits a prescribed-time ISS Lyapunov function $V_{i}(\chi_{i})$
in the form of Definition \ref{defn:PTISS} with $d_{i}$ and $\chi_{3-i}$
as the inputs, i.e.,
\begin{eqnarray}
&V_{1}(\chi_{1})\sim\left\{ \underline{\alpha}_{1}, \overline{\alpha}_{1},\alpha_{1},[V_{2}(\chi_{2}),l_{1}\alpha_{2}],[\epsilon_{1}(\|d_{1}\|),\gamma_{1}]\right\} \label{eq:PTISS-1}\\
&V_{2}(\chi_{2})\sim\left\{ \underline{\alpha}_{2}, \overline{\alpha}_{2},\alpha_{2},[V_{1}(\chi_{1}),l_{2}\alpha_{1}],[\epsilon_{2}(\|d_{2}\|),\gamma_{2}]\right\} \label{eq:PTISS-2}
\end{eqnarray}
 where $\underline{\alpha}_{i}$, $ \overline{\alpha}_{i}$, $\alpha_{i}$,
$\gamma_{i}$, $\epsilon_{i}\in\mathcal{K}_{\infty}$, and $l_{i}$
is positive a finite constant for $i=1,2$. If
\begin{gather}
l_{1}l_{2}<1,\label{eq:k}\\
\min\{\alpha_{1}(s),\alpha_{2}(s)\}=\mathcal{O}\left[\gamma_{1}(s)+\gamma_{2}(s)\right]\label{eq:O_s}
\end{gather}
 hold for $s\in\mathbb{R}_{+}$, then $\chi$ is uniformly bounded
for $t\in\mathcal{T}_{p}$ with upper bound
\begin{equation}
\|\chi\|\leq\tilde{\alpha}(\|\chi_{0}\|)\label{eq:bound_chi}
\end{equation}
 where $\tilde{\alpha}$ is some $\mathcal{K}_{\infty}^{e}$ function.
\ethm

\brem
In contrast to asymptotic stabilization ~\citep{jiang2018small31}
and finite-time stabilization ~\citep{pavlichkov2018finite29,li2023lyapunov30},
the prescribed-time stabilization of an interconnected system
requires not only that the coupling parts of prescribed-time ISS Lyapunov functions
meet the small-gain condition (\ref{eq:k}), but also
specific condition on the relationship between prescribed-time convergent gain and prescribed-time ISS gain with respect to external disturbance, see (\ref{eq:O_s}) for details.
\erem

\brem For Theorem \ref{lem:SG-Criterion}, if (\ref{eq:k}) holds and
(\ref{eq:O_s}) is modified to
$\min\{\alpha_{1}(s),\alpha_{2}(s)\}=\mathcal{O}\left[\tilde{\gamma}(s)(\gamma_{1}(s)+\gamma_{2}(s))\right]$
for some $\tilde{\gamma}\in\mathcal{K}_{\infty}$, then all results
in Theorem \ref{lem:SG-Criterion} still hold. Additionally, $\chi$  satisfies
$\|\chi(t)\|\leq(\tilde{\gamma}'\circ\tilde{\gamma}(\mu))^{-1}\tilde{\alpha}'(\|\chi(t_{0}))\|$
for some $\tilde{\gamma}'\in\mathcal{K}_{\infty}$ and $\tilde{\alpha}'\in\mathcal{K}_{\infty}^{e}$,
which implies that $\chi$ achieves prescribed-time convergence towards
zero  in the presence of any bounded external disturbances.
This highlights
the robustness advantage of prescribed-time control over asymptotic and finite-time control.  \erem

\section{DPTCO Implementation \label{sec:DPTCO-Implementation}}
Through the application of time-varying state transformation \eqref{eq:tilde_e_r} and Lemma \ref{lem:PT-Mapping}, the implementation of the DPTCO, as well as the uniform boundedness of all internal signals, now hinges on the boundedness of $\tilde{e}_{r}$
and $\tilde{e}_{s}$. In this section, we elaborate on the parameter design and provide a detailed stability analysis of the dynamics of $\tilde{e}_{r}$
and $\tilde{e}_{s}$.

\subsection{The Prescribed-Time ISS Lyapunov Functions}
Due to the introduction of adaptive estimation $\hat \theta$, it is essential to incorporate a positive definite, differentiable term with respect to estimation error $\tilde \theta$ in the corresponding Lyapunov function.
Therefore, we define
\begin{equation}
\tilde{e}_{s}'=[\tilde{e}_{s}\t,\tilde{\theta}\t]\t.\label{eq:tilde-e-s'}
\end{equation}
Then according to Lemma \ref{lem:PT-Mapping} and Theorem \ref{lem:SG-Criterion},
for $t\in\mathcal{T}_{p}$, the DPTCO is achieved if the dynamics of $\tilde{e}_{r}$-
and $\tilde{e}_{s}'$ admit prescribed-time ISS Lyapunov
functions as described in (\ref{eq:PTISS-1}) and (\ref{eq:PTISS-2}). Additionally,  the design parameter $\eta$ in (\ref{eq:eta}) must ensure that conditions (\ref{eq:k})
and (\ref{eq:O_s}) are satisfied. For simplicity, define $\bar{k}_{ \overline{m}}=\max\{k_{ \overline{m}  }^{1},\cdots,k_{  \overline{m}  }^{N}\}$,  $\underline{k}_{\underline{m}}=\min\{k_{\underline{m}}^{1},\cdots,k_{\underline{m}}^{N}\}$ and
\begin{eqnarray}
\begin{gathered}\delta=\max\{4/\lambda_{2},(4\varrho_{c}^{2}+1)/\rho_{c}\},\\
\bar{\varepsilon}=\max\left\{ 1,\bar{k}_{ \overline{m}  }\right\} ,\quad\underline{\varepsilon}=\min\left\{ 1,\underline{k}_{\underline{m}}\right\},\\
\bar{\delta}=\delta\max\{1,\lambda_{2}^{-1}\}/2+1/2,\quad \underline{\delta}=\delta\min\{1,\lambda_{N}^{-1}\}/2.
\end{gathered}
\label{eq:delta}
\end{eqnarray}
Then we proposed the following design criteria for parameter $\eta$
in (\ref{eq:eta}).

$\textbf{DC}_{\textbf{1}}$: $\iota$ in \eqref{eq:z_1^i} is any any constant satisfying
$\iota>2$;

$\textbf{DC}_{\textbf{2}}$: $c$ in \eqref{eq:dot_varpi_i} and \eqref{eq:dot_v_i} is designed as
\begin{equation}
c=2\bar{\delta}c^{*}+4\iota\tilde{b}+4\bar{\delta}\iota\tilde{b}\label{eq:c}
\end{equation}
 with any $c^{*}>0$, where $\tilde{b}$ and $\bar{\delta}$ are denoted
in (\ref{eq:bound_mu_derivative}) and (\ref{eq:delta}), respectively;

$\textbf{DC}_{\textbf{3}}$: For $i\in\mathcal{V}$, $\sigma^{i}$ in \eqref{eq:dot_hat_theta^i}
is designed to such that
\begin{equation}
\sigma_{\min}>\max\{\tilde{b},c_{\Delta}c_{s}/(c^{*}\underline{\delta})\}\label{eq:sigma}
\end{equation}
 where $\sigma_{\min}$ has been denoted in \eqref{eq:sigma-min}, and
 \begin{equation}
 \begin{aligned}
 c_{\Delta}&=c\big(\varrho_{c}^{2}\delta/(2\rho_{c}) +\lambda_{N}^{2}\delta/(2\lambda_{2})+2\delta^{2}\\
 &\quad +\varrho_{c}^{2}/(2\rho_{c})+2\varrho_{c}^{2}\big), \\
 c_{s}&=4c^{2}(2\lambda_{N}+\varrho_{c}+1)^{2}
 \end{aligned}\label{eq:c-delta}
 \end{equation}

$\textbf{DC}_{\textbf{4}}$: For $i\in\mathcal{V}$, $k_{1}^{i}$
and $k_{2}^{i}$ in \eqref{eq:pi^i} are designed as
\begin{equation}
k_{1}^{i}=k_{1}^{i,*}+\tilde{b}\iota+1,\quad k_{2}=k_{2}^{i,*}+\frac{1}{2}(k_{1}^{i})^{2} (k_{  \overline{m} }^{i})^{2}+\frac{1}{2}\label{eq:k-1}
\end{equation}
 with $k_{1}^{i,*}$ and $k_{2}^{i,*}$ are chosen to satisfy
\begin{equation}
\underline{k}_{1}^{*}>c_{s}(c_{\Delta}/(c^{*}\underline{\delta})+1)/2,\quad k_{2}^{i,*}>c_{\Delta}c_{s}k_{ \overline{m}  }^{i}/(2c^{*}\underline{\delta})\label{eq:k_12^*}
\end{equation}
 where $\underline{k}_{1}^{*}=\min\{k_{1}^{i,*},\cdots,k_{1}^{N,*}\}$.

In the following,
we show that there exist prescribed-time ISS Lyapunov
functions for $\tilde{e}_{r}$- and $\tilde{e}_{s}'$-dynamics.

\blem\label{the:1} Suppose Property \ref{prop:1-1}-\ref{prop:1-3},
Assumption \ref{ass:graph}, \ref{ass:solvable} hold and the design parameter
$\eta$ satisfies $\textbf{DC}_{\textbf{1}}$-$\textbf{DC}_{\textbf{4}}$,
let $\tilde{e}_{r}=\chi_{1}$, $\tilde{e}_{s}'=\chi_{2}$, $d_{1}=0$
and $d_{2}=\theta$, then $\tilde{e}_{r}$-dynamics in (\ref{eq:dot_til_e_r})
and $\tilde{e}_{s}'$-dynamics in (\ref{eq:dot_til_e_s}), (\ref{eq:dot_til_theta})
admit the prescribed-time ISS Lyapunov functions as described in (\ref{eq:PTISS-1}) and
(\ref{eq:PTISS-2}) with
\begin{eqnarray}
 && \begin{aligned} & \underline{\alpha}_{1}(s)=\underline{\delta}s^{2},\;\; \overline{\alpha}_{1}(s)=\bar{\delta}s^{2},\;\;\alpha_{1}(s)=c^{*}s,\\
 & l_{1}=c_{\Delta}/\tilde{k},\;\;\gamma_{1}(s)=s,\;\;\epsilon_{1}(\|d_{1}\|)=0,
\end{aligned}
\label{eq:con1}\\
 && \begin{aligned} & \underline{\alpha}_{2}(s)=\underline{\varepsilon}s^{2},\;\; \overline{\alpha}_{2}(s)=\bar{\varepsilon}s^{2},\;\;\alpha_{2}(s)=\tilde{k}s,\\
 & l_{2}=c_{s}/(c^{*}\underline{\delta}),\;\; \gamma_{2}(s)=s,\;\;\epsilon_{2}(\|d_{2}\|)=\sigma_{\max}\|d_{2}\|^{2}
\end{aligned}
\label{eq:con2}
\end{eqnarray}
 where
$
\tilde{k}=\min\{\tilde{k}_{1},\tilde{k}_{2},\sigma_{\min}\}%%\label{eq:k*}
$
with $\tilde{k}_{1}=2\underline{k}_{1}^{*}-c_{s}$, $\tilde{k}_{2}=2\min\{k_{2}^{1,*}/k_{\overline{m}}^{1},\cdots,k_{2}^{N,*}/k_{\overline{m}}^{N}\}$.\elem

The proof of Lemma \ref{the:1} is given in appendix. 

\subsection{Prescribed-Time Stability of Closed-Loop Systems}
The following theorem is to prove the DPTCO objectives in (\ref{eq:objective1})
and (\ref{eq:objective2}) are achieved by the proposed algorithm.

\bthm\label{thm-1}
Consider the closed-loop systems consist of (\ref{eq:Euler-Lag}),
(\ref{eq:dot_varpi_i}), (\ref{eq:dot_v_i}), (\ref{eq:dot_varpi_v}),
(\ref{eq:tau^i}) and (\ref{eq:dot_hat_theta^i}) for $i\in\mathcal{V}$ under Property \ref{prop:1-1} - \ref{prop:1-3} and Assumption \ref{ass:graph} - \ref{ass:cost_func} .
Suppose the parameter $\eta$ is selected to satisfy $\textbf{DC}_{\textbf{1}}$ - $\textbf{DC}_{\textbf{4}}$,
then the DPTCO problem is solved in the sense that the output $y^{i}$
of each subsystem achieves prescribed-time convergence towards the optimal
solution of the global objective function (\ref{eq:opti_problem}) at $T+t_{0}$
and remains at the optimal solution afterwards. Moreover, all internal
signals in the closed-loop systems are uniformly bounded. \ethm

\proofnow Choosing any $\iota>2$ and $c^{*}>0$, we have  that $c$ in
(\ref{eq:c}) is fixed, and then $l_{2}$ in (\ref{eq:con2}) is fixed.
For $\sigma^{i}$, $k_{1}^{i}$ and $k_{2}^{i}$ satisfy $\textbf{DC}_{\textbf{3}}$
and $\textbf{DC}_{\textbf{4}}$, by some algebra operations, one has
$\tilde{k}>c_{\Delta}c_{s}/(c^{*}\underline{\delta})$, which further
leads to
\[
c_{\Delta}c_{s}/(\tilde{k}c^{*}\underline{\delta})<1.
\]
Thus (\ref{eq:k}) is achieved according to (\ref{eq:con1}) and (\ref{eq:con2}).
And we have
\begin{align*}
 & \sup_{s\in\mathbb{R}_{+}}\left[(\gamma_{1}(s)+\gamma_{2}(s))/\min\{\alpha_{1}(s),\alpha_{2}(s)\}\right]\\
 & =\sup_{s\in\mathbb{R}_{+}}\left[2s/\min\{c^{*}s,\tilde{k}s\}\right]=2/\min\{c^{*},\tilde{k}\}<\infty
\end{align*}
 which means (\ref{eq:O_s}) is achieved.
 Define $\eta_1 = (c_\Delta c_s/(\tilde k c^* \underline \delta))^{1/2}$, $\eta_2 = 2\eta_1 \sigma_{\max} /(\min\{c^*,\tilde k\}(1-\eta_1))$, and $\eta_3 = \max\{c_s/(c^*\underline \delta),\eta_1\}(\overline \delta +\overline \varepsilon)$.
 According to Theorem \ref{lem:SG-Criterion} and Lemma \ref{the:1}, we have
\begin{equation}
\|[\tilde{e}_{r}\t,(\tilde{e}_{s}')\t]\|\leq\tilde{\alpha}(\|[\tilde{e}_{r}\t(t_{0}),(\tilde{e}_{s}'(t_{0}))\t]\|)\label{eq:bound_tilde_e_r_s}
\end{equation}
 for $\mathcal{K}_{\infty}^{e}$ function $\tilde{\alpha}$ denoted
as
\begin{align}
\tilde{\alpha}(s)=\big(\underline \delta^{-1/2}+\underline \varepsilon^{-1/2}\big)\big(\eta_2^{1/2}\|\theta\| +\eta_3^{1/2}s  \big).\label{eq:tilde-alpha}
\end{align}
By Lemma \ref{lem:PT-Mapping}, the DPTCO is achieved and all internal
signals in the closed-loop systems are uniformly bounded for $t\in\mathcal{T}_{p}$.

Note that for $t\in\tilde{\mathcal{T}}_{p}$, by (\ref{eq:dot_varpi_v})
\begin{align*}
\left[\begin{array}{c}
\varpi(t)\\
v(t)
\end{array}\right] & =\left[\begin{array}{c}
\varpi(T+t_{0})\\
v(T+t_{0})
\end{array}\right]+\int_{T+t_{0}}^{t}\left[\dot{\varpi}\t(\tau),\dot{v}\t(\tau)\right]\t\mathrm{d}\tau\\
 & =\lim_{t\to T+t_{0}}\left[\begin{array}{c}
\varpi(t)\\
v(t)
\end{array}\right]=\left[\begin{array}{c}
1_{N}\otimes z^{*}\\
-\nabla F(1_{N}\otimes z^{*})
\end{array}\right]
\end{align*}
By (\ref{eq:tau^i}) and (\ref{eq:bound_pi}),
one has
\[
\lim_{t\to(T+t_{0})^{-}}\tau(t)=0=\tau(T+t_{0})=\lim_{t\to(T+t_{0})^{+}}\tau(t)
\]
 which means controller $\tau(t)$ is continuous at $t=T+t_{0}$,
and further, is continuous on the global time interval. By (\ref{eq:Euler-Lag})
and (\ref{eq:tau^i}), for $t\in\mathcal{T}_{p}$, one has $\ddot{q}=-M^{-1}(q)C(q,\dot{q})\dot{q}$.
Solving the differential equation yields
\begin{align*}
\dot{q}(t) & =\exp\left(-\sint_{T+t_{0}}^{t}M^{-1}(q(\tau))C(q(\tau),\dot{q}(\tau))\mathrm{d}\tau\right)\dot{q}(T+t_{0})\\
 & =0,\quad\forall t\in\tilde{\mathcal{T}}_{p}
\end{align*}
where we used $\dot{q}(T+t_{0})=\lim_{t\to T+t_{0}}\dot{q}(t)=0$.
Then for $t\in\tilde{\mathcal{T}}_{p}$, $\dot{y}=\dot{q}=0$ and
\[
y(t)=y(T+t_{0})+\sint_{T+t_{0}}^{t}\dot{y}(\tau)\mathrm{d}\tau=1_{N}\otimes z^{*}
\]
where we used $y(T+t_{0})=\lim_{t\to T+t_{0}}y(t)=1_{N}\otimes z^{*}$.
The above equations mean $e_{r}(t)=0$ and $e_{s}(t)=0$ hold for
$t\in\tilde{\mathcal{T}}_{p}$.

Recall the definition of $\tilde e_s'$ in \eqref{eq:tilde-e-s'}, one has $\|[\tilde{e}_{r}\t,\tilde{e}_{s}\t]\|\leq \|[\tilde{e}_{r}\t,(\tilde e_s')\t\|$.
By (\ref{eq:bound_tilde_e_r_s}) and \eqref{eq:tilde-alpha}
\begin{align}
\|[\tilde{e}_{r}\t,\tilde{e}_{s}\t]\|&\leq \tilde{\alpha}(\|[\tilde{e}_{r}\t(t_{0}),(\tilde{e}_{s}'(t_{0}))\t]\|)\notag \\
&\leq \tilde{\alpha}'(\|e_{r}\t(t_{0}),e_{s}\t(t_{0})\|)\notag
\end{align}
with $\tilde{\alpha}'$ belongs to $\mathcal{K}_{\infty}^{e}$ denoted
as $\tilde{\alpha}'(s)=\big(\underline \delta^{-1/2}+\underline \varepsilon^{-1/2}\big)\big(\eta_2^{1/2}\|\theta\| +\eta_3^{1/2} \|\tilde \theta(t_0)\|+\eta_3^{1/2}s  \big)$
Then by (\ref{eq:bounds_e_r_s}) in Lemma \ref{lem:PT-Mapping}, we
have
\[
\|e_{r}\t(t),e_{s}\t(t)\|\leq\beta_{e}(\|[e_{r}\t(t_{0}),e_{s}\t(t_{0})]\|,t-t_{0})
\]
 with $\beta_{e}(r,s)\in\mathcal{KL}_{T}^{e}$ denoted as
\[
\beta_{r}(r,s)=\left\{ \begin{aligned} & \mu^{-\iota+1}(s+t_{0})\left[\mu^{-1}(t_{0})\gamma_{r}\circ\tilde{\alpha}'(r)\right.\\
 & \left.+\gamma_{s}\circ\tilde{\alpha}'(r)\right],\quad0\leq r,\;0\leq s<T,\\
 & 0,\quad0\leq r,\;T\leq s,
\end{aligned}
\right.
\]
 which implies the prescribed-time stability of $[\tilde{e}_{r}\t,\tilde{e}_{s}\t]$-dynamics
according to Definition \ref{def:PT-Stbility}. This completes the proof. \eproof

\section{Simulation Results \label{sec:Simulation}}

\begin{figure}
\includegraphics[scale=0.15]{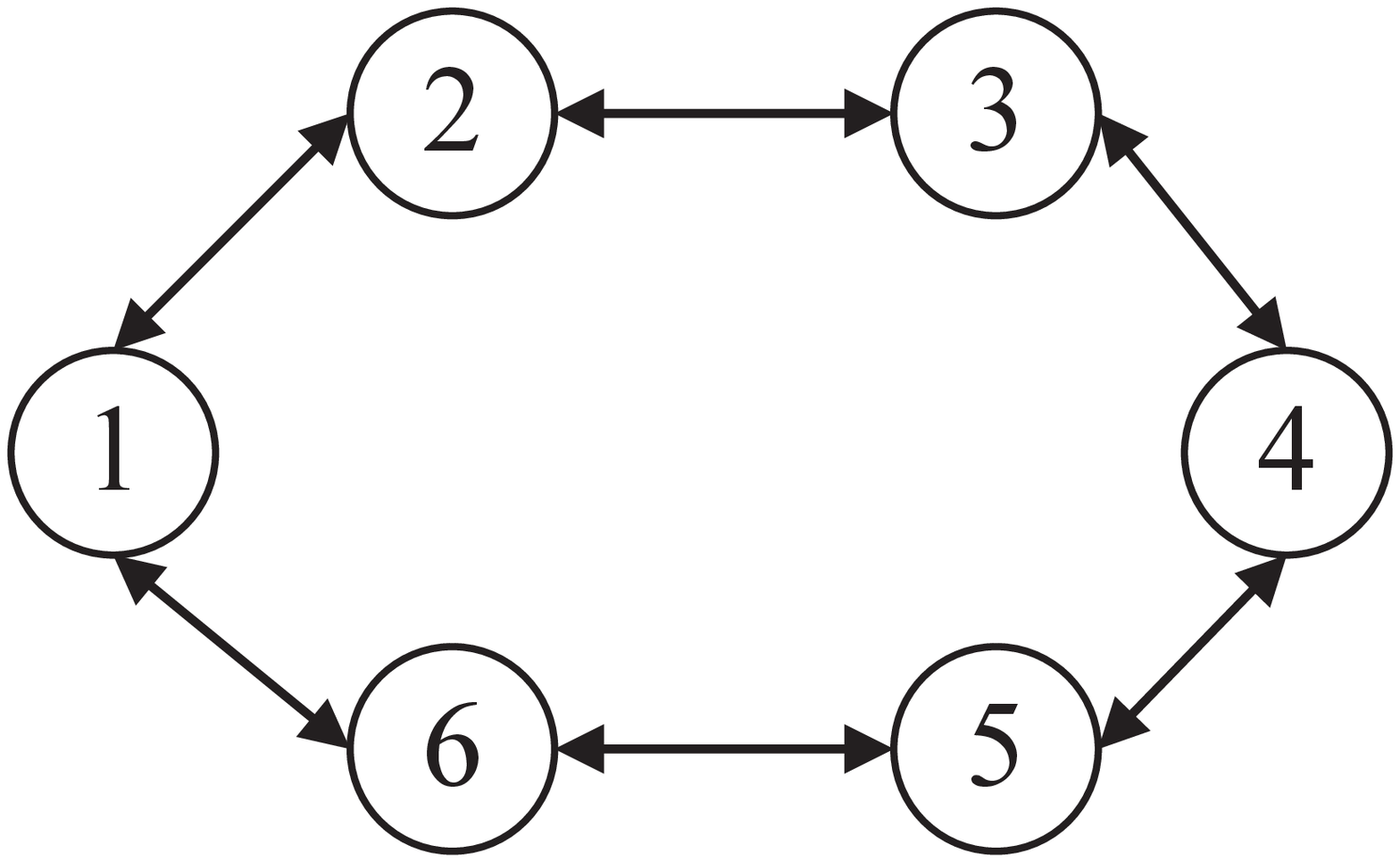}

\caption{\label{fig:graph} Communication graph $\mathcal{G}$}
\end{figure}

The information flow among all agents is described by a undirected graph in Fig.
\ref{fig:graph}. Each subsystem in (\ref{eq:Euler-Lag}) is in mathematical
model of
a two-link
revolute joint robotic manipulator given in \cite{zhang2017distributed05}
with $q^{i}=[q_{1}^{i},q_{2}^{i}]\t$ and
\begin{gather*}
M^{i}(q^{i})=\left[\begin{array}{cc}
p_{1}^{i}+2p_{2}^{i}\cos q_{2}^{i} & p_{3}^{i}+p_{2}^{i}\cos q_{2}^{i}\\
p_{3}^{i}+p_{2}^{i}\cos q_{2}^{i} & p_{3}^{i}
\end{array}\right],\\
C^{i}(q^{i},\dot{q}^{i})=\left[\begin{array}{cc}
-p_{2}^{i}(\sin q_{2}^{i})\dot{q}_{2}^{i} & -p_{2}^{i}(\sin q_{2}^{i})(\dot{q}_{1}^{i}+\dot{q}_{2}^{i})\\
p_{2}^{i}(\sin q_{2}^{i})\dot{q}_{1} & 0
\end{array}\right].
\end{gather*}
Although we do not know the exact value of $M^{i}(q^{i})$ and $C^{i}(q^{i},\dot{q}^{i})$
for $i=6$, it is assumed that the unknown parameters are in the following
ranges $p_{1}^{i}\in[1,2]{\rm kg\cdot m^{2}}$, $p_{2}^{i}\in[0.05,0.2]{\rm kg\cdot m^{2}}$
and $p_{3}^{i}\in[0.2,0.5]{\rm kg\cdot m^{2}}$. Obviously, $\dot{M}^{i}(q^{i})-2C(q,\dot{q})$
is  skew-symmetric. The matrix $\Omega^{i}(q^{i},\dot{q}^{i},x,y)$
and $\theta^{i}$ in Property \ref{prop:1-2} can be calculated as
$\Omega^{i}(q^{i},\dot{q}^{i},x,y)=[x_{1},\cos q_{2}^{i}(2x_{1}+x_{2})-\sin q_{2}^{i}(\dot{q}_{2}^{i}y_{1}+\dot{q}_{1}^{i}y_{2}+\dot{q}_{2}^{i}y_{2}),x_{2};0,(\cos q_{2}^{i})x_{1}+(\sin q_{2}^{i})\dot{q}_{1}^{i}y_{1},x_{1}+x_{2}]$
and $\theta^{i}=[p_{1}^{i},p_{2}^{i},p_{3}^{i}]\t$ with $x=[x_{1},x_{2}]\t$,
$y=[y_{1},y_{2}]\t$. Simple calculation verifies that Properties
\ref{prop:1-2} and \ref{prop:1-3} is satisfied with $\rho^{i}=4.24$,
$k_{\bar{m}}^{i}=2.63$ and $k_{\underline{m}}^{i}=0.0789$.

\begin{figure}
\includegraphics[scale=0.6]{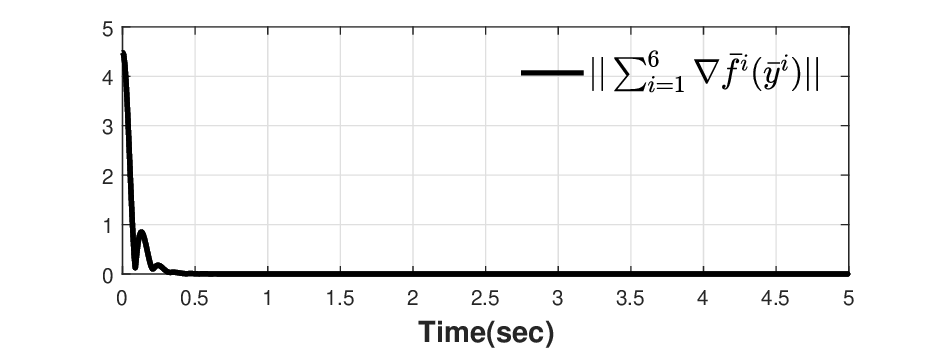}

\caption{\label{fig:gradient} Trajectory of the norm of global objective function,
where $\nabla\bar{f}^{i}(\bar{y}^{i})$ denotes the gradient of local
objective function in (\ref{eq:opti-prob}).}
\end{figure}

\begin{figure}

\includegraphics[scale=0.6]{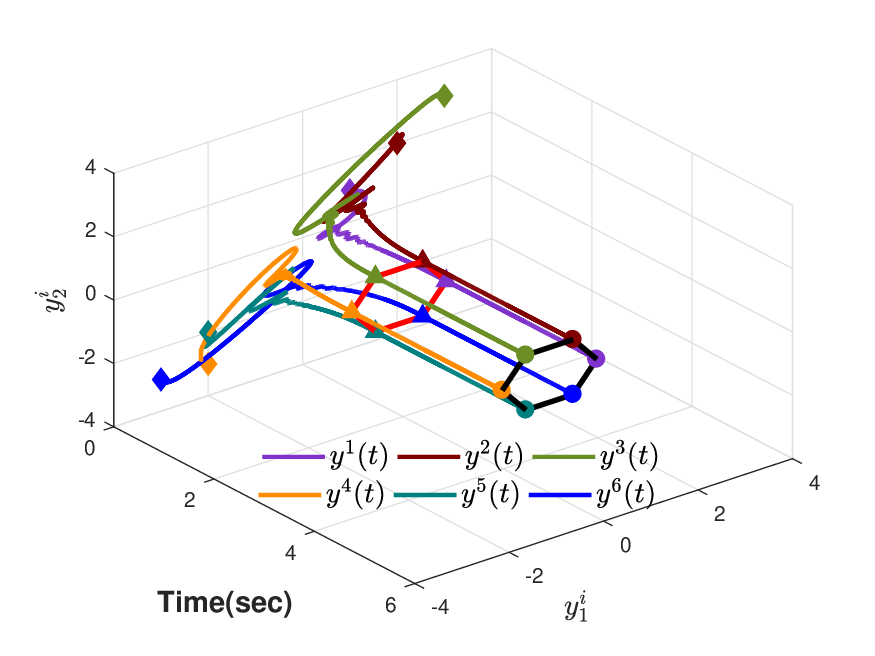}\caption{\label{fig:Formation} Trajectories of positions $y^{i}$ of the six
robots in three-dimensional space for $0\protect\leq t\protect\leq5s$,
where $\blacklozenge$, $\blacktriangle$ and $\bullet$ denote the
initial position, position at the prescribed time, and final position,
respectively, the red and black regular hexagons represent the formation
of the six robots at the prescribed time and final time.  }

\end{figure}

We consider a network composed of six EL agents to cooperatively search for
an unknown heat source. Each agent is equipped with sensors to measure
gradient information of heat source with respect to distance. The
objective is to design controller $u^{i}$ such that the six robots
approach the heat source in a formation, and reduce the total displacement
of the six robots from their original location. Thus, the global objective function is designed as
\begin{equation}
\begin{aligned} & \min\ssum_{i=1}^{6}\left(s_{1}^{i}\|y^{i}-d^{*}\|^{2}+s_{2}^{i}\|y^{i}-y^{i}(t_{0})\|^{2}\right),\\
 & \mbox{subject to}\quad y^{i}-y^{j}=\omega^{i}-\omega^{j}
\end{aligned}
\label{eq:opti_problem-1}
\end{equation}
where $d^{*}$ denotes the two-dimensional coordinates of the heat
source, $\omega^{1}=[1,0]\t$, $\omega^{2}=[1/2,\sqrt{3}/2]\t$, $\omega^{3}=[-1/2,\sqrt{3}/2]\t$,
$\omega^{4}=[-1,0]\t$, $\omega^{5}=[-1/2,-\sqrt{3}/2]\t$, $\omega^{6}=[1/2,-\sqrt{3}/2]\t$
represent the formation shape, and $s_{1}^{i}$ and $s_{2}^{i}$ are
objective weights.

By defining $\bar{y}^{i}=y^{i}-\omega^{i}$, the optimization problem
(\ref{eq:opti_problem-1}) is transformed into
\begin{equation}
\begin{aligned} & \min\ssum_{i=1}^{6}\left(s_{1}^{i}\|\bar{y}^{i}+\omega^{i}-d^{*}\|^{2}+s_{2}^{i}\|\bar{y}^{i}+\omega^{i}-y^{i}(t_{0})\|^{2}\right),\\
 & \mbox{subject to}\quad\bar{y}^{i}=\bar{y}^{j}
\end{aligned}
\label{eq:opti-prob}
\end{equation}
which is consistent with (\ref{eq:opti_problem}). For optimization
problem (\ref{eq:opti-prob}), we design $\varpi^{i}$- and $v^{i}$-dynamics
as in the form of (\ref{eq:dot_varpi_i}), (\ref{eq:dot_v_i}) and
(\ref{eq:dot_varpi_v}) such that $\varpi^{i}$ converges to optimal
solution with the prescribed time and remains at the
optimal solution thereafter.
Then, the reference trajectory of each robot's
dynamics is changed as
\[
\varpi^{i,'}=\varpi^{i}+\omega^{i}.
\]
 Replacing $\varpi^{i}$ in (\ref{eq:tau^i})-(\ref{eq:varsigma^i})
with $\varpi^{i,'}$, the proposed algorithm can solve optimization
problem (\ref{eq:opti_problem-1}). Let the initial condition be
$y^{1}(t_{0})=[1,1]\t$, $y^{2}(t_{0})=[2,2]\t$, $y^{3}(t_{0})=[3,3]\t$,
$y^{4}(t_{0})=[-2,-3]\t$, $y^{5}(t_{0})=[-2,-2]\t$, $y^{6}(t_{0})=[-3,-3]$,
$\dot{q}^{i}(t_{0})=[0,0]\t$, $\varpi^{i}(t_{0})=y^{i}(t_{0})$,
$v^{i}(t_{0})=[0,0]\t$, $\hat{\theta}^{i}(t_{0})=[2,2,2]\t$ for
$i=1,\cdots,6$. The initial time $t_{0}$ is set
to $t_{0}=0$, the
prescribed-time scale $T=2s$ and total simulation time is $5s$.
The parameters and gain function are chosen as $c=1.3$, $\iota=2.44$,
$k_{1}^{i}=5$, $k_{2}^{i}=30$, $\sigma^{i}=4$ for $i=1,\cdots,6$,
$\mu=10/(T+t_{0}-t)$. The weight coefficients $s_{1}^{i}$ and $s_{2}^{i}$
for objective function (\ref{eq:opti_problem-1}) are chosen as $s_{1}^{i}=1$
and $s_{2}^{i}=1$ for $i=1,\cdots,6$. The coordinate of the heat source
is set as $d^{*}=[0,0]\t$. The real value of unknown parameters in
$M^{i}(q^{i})$ are $p_{1}^{i}=1.301{\rm kg\cdot m^{2}}$, $p_{2}^{i}=0.056{\rm kg\cdot m^{2}}$
and $p_{3}^{i}=0.296{\rm kg\cdot m^{2}}$ for $i=1,\cdots,6$. The
fourth-order Runge-Kutta method is used to solve the ordinary differential
equations, and the sampling time is $10^{-5}s$.

The simulation results are shown in Fig. \ref{fig:gradient} and Fig.~\ref{fig:Formation}. Fig.~\ref{fig:gradient}
demonstrates that
$\bar{y}^{i}$ converges to the optimal solution of the optimization problem (\ref{eq:opti-prob})
within a prescribed time and
remains at the
optimal solution
% afterwards
for $i=1,\cdots,6$, which further implies the optimization problem (\ref{eq:opti_problem-1}) is solved within
% a prescribed time.
the given time frame.
In Fig.~\ref{fig:Formation}, the six robots approach the heat source in formation within a prescribed time and maintain formation subsequently.
% afterwards.

\section{Conclusion\label{sec:Conclusion}}

In this paper, we propose a novel DPTCO algorithm for a class of NELSs.
The position of each robot in the network achieves prescribed-time
convergence towards the optimal solution within a prescribed time and
remains at the
optimal solution afterwards. The optimal solution search
is in a
feedback loop
in which position-dependent measured gradient
value is used to construct the auxiliary systems. A prescribed-time
small-gain criterion is proposed to characterize the interconnected
error system and is also used to guide the selection of control parameters.
As a future work, it
would be
intriguing
to further consider the DPTCO where the
local objective functions subject to bound, equality, and inequality constraints.

\scriptsize{
\bibliographystyle{apalike}        % Include this if you use bibtex
\bibliography{ref}       }

\appendix
\normalsize{
\section{Proof of Lemma \ref{lem:PT-Mapping}}
Due to (\ref{eq:tilde_e_r}), $e_{r}$ satisfies
\begin{equation}
\|e_{r}\|\leq\mu^{-\iota}\|\tilde{e}_{r}\|\leq\mu^{-\iota}\|\tilde{e}_{r}\|_{\mathcal{T}}.\label{eq:norm_e_r}
\end{equation}
 Thus the first inequality in (\ref{eq:bounds_e_r_s}) holds with
$\gamma_{r}(s)=s$. The $\tilde{e}_{s}$ in (\ref{eq:tilde_e_r})
can be rewritten in a detailed form as
\begin{equation}
\tilde{e}_{s}=\left[\begin{array}{c}
\mu^{\iota}e_{y}\\
\mu^{\iota-1}\left(\dot{q}+\mu\bar{K}_{1}e_{y}\right)
\end{array}\right].\label{eq:tilde_e_s-1}
\end{equation}
Then one has
\begin{equation}
\|e_{y}\|\leq\mu^{-\iota}\|\tilde{e}_{s}\|_{\mathcal{T}}\label{eq:norm_e_y}
\end{equation}
where we used $\mu^{\iota}\|e_{y}\|\leq\|\tilde{e}_{s}\|\leq\|\tilde{e}_{s}\|_{\mathcal{T}}$.
According to (\ref{eq:tilde_e_s-1}) and (\ref{eq:norm_e_y}), $\dot{q}$
satisfies
\begin{align}
\|\dot{q}\|
 & \leq\mu^{1-\iota}\|\tilde{e}_{s}\|+\mu\|\bar{K}_{1}\|\|e_{y}\|\nonumber \\
 & \leq\mu^{1-\iota}\|\tilde{e}_{s}\|_{\mathcal{T}}+\|\bar{K}_{1}\|\mu^{1-\iota}\|\tilde{e}_{s}\|_{\mathcal{T}}\nonumber \\
 & \leq(1+\|\bar{K}_{1}\|)\mu^{1-\iota}\|\tilde{e}_{s}\|_{\mathcal{T}}.\label{eq:norm_dot_q}
\end{align}
By (\ref{eq:norm_e_y}) and (\ref{eq:norm_dot_q}), the second inequality
in (\ref{eq:bounds_e_r_s}) holds with $\gamma_{s}(s)=(1+\|\bar{K}_{1}\|+b^{-1})s$ where we used $\mu^{-1}(t)\leq \mu^{-1}(t_{0})=b^{-1}$.
By (\ref{eq:bound_dot_e_r}) and \eqref{eq:bounds_e_r_s},
the upper bound of $\dot e_r$ is
\begin{equation}
\|\dot{e}_{r}\|\leq\mu^{-\iota+2}\gamma_{e}(\|\tilde{e}\|_{\mathcal{T}})\label{eq:bound_dot_e_r-1}
\end{equation}
 where $\tilde{e}=\left[\tilde{e}_{r}\t,\tilde{e}_{s}\right]\t$ and
$\gamma_{e}$ is a $\mathcal{K}_{\infty}$ function denoted as $\gamma_{e}(s)=c(2\lambda_{N}+\varrho_{c}+1)(b^{-1}\gamma_r (s)+ \gamma_s(s))$.
Since $\iota>2$, (\ref{eq:bound_dot_e_r-1}) means that $\dot{\varpi}^{i}$
and $\dot{v}^{i}$ are uniformly bounded if $\tilde{e}_{r}$ and $\tilde{e}_{s}$
are uniformly bounded.

Define
$z_{j}=[(z_{j}^{1})\t,\cdots,(z_{j}^{N})\t]\t$, where
$z_{j}^{i}$ is given by  (\ref{eq:z_1^i}) and  (\ref{eq:z_2^i}) for $j=1,2$ and $i=1,\cdots,N$.
Based on (\ref{eq:bound_mu_derivative}), (\ref{eq:norm_e_y}) and (\ref{eq:norm_dot_q}),
the upper bound of $z_{j}$
\[
\|z_{j}\|\leq\mu^{2-\iota}\gamma_{z,j}(\|\tilde{e}_{s}\|_{\mathcal{T}}), \quad j=1,2
\]
where $\gamma_{z,1}(s)=(\tilde{b}(\iota-1)+\|\bar{K}_{1}\|)(1+\|\bar{K}_{1}\|)s+\iota\tilde{b}\|\bar{K}_{1}\|s$,
$\gamma_{z,2}(s)=\|\bar{K}_{1}\|b^{-1}s$. By Property \ref{prop:1-2}, let $\omega^i = [(\dot q^i)\t,(z_1^i)\t, (z_2^i)\t]\t$,
$\Omega(q,\dot{q},z_{1},z_{2})$
satisfies
\begin{equation}
\|\Omega(q,\dot{q},z_{1},z_{2})\|\leq\mu^{-\iota+2}\gamma_{\Omega}(\|\tilde{e}_{s}\|_{\mathcal{T}})\label{eq:bound_Omega}
\end{equation}
with
\begin{align*}
\gamma_{\Omega}(s) & =N\bar \rho \big(b^{-1} \gamma_s (s)+b^{-\iota} (\gamma_s(s))^2 + \ssum_{j=1}^2 \gamma_{z,j}(s)\big) \\
&\quad + N\bar \rho \ssum_{j=1}^2 b^{-\iota+2}(\gamma_{z,j}(s))^2
\end{align*}
where $\bar{\rho}=\max\{\rho^{1},\cdots,\rho^{N}\}$ with $\rho^{i}$
been denoted in Property \ref{prop:1-2} and we used $\left(\ssum_{i=1}^{n}|x^{i}|\right)^{p}\leq\ssum_{i=1}^{n}\left(|x^{i}|\right)^{p}$
for $x^{i}\in\mathbb{R}$, $0<p\leq1$ and $i=1,\cdots,n$. By (\ref{eq:tilde_e_s-1})
and (\ref{eq:bound_Omega}), the upper bound of $\varsigma=\left[(\varsigma^{1})\t,\cdots,(\varsigma^{N})\t\right]\t$
with $\varsigma^{i}$ denoted in (\ref{eq:varsigma^i}) is

\begin{equation}
\|\varsigma(q,\dot{q},\varpi,\mu)\|\leq\mu\gamma_{\varsigma}(\|\tilde{e}_{s}\|_{\mathcal{T}})\label{eq:bound_varsigma}
\end{equation}
 where $\gamma_{\varsigma}(s)=2s\gamma_{\Omega}(s)$. Define the Lyapunov
function candidate for $\hat{\theta}$-dynamics as $U(\hat{\theta})=\hat{\theta}\t\hat{\theta}/2$,
 then its time derivative along trajectory of (\ref{eq:dot_hat_theta^i})
is
\begin{align}
\dot{U}(\hat{\theta}) & =-\mu\hat{\theta}\t\bar \sigma\hat{\theta}+\hat{\theta}\t\varsigma(q,\dot{q},\varpi,\mu)\nonumber \\
 & \leq-\sigma_{\min}\mu U(\hat{\theta})+(2\sigma_{\min})^{-1}\mu(\gamma_{\varsigma}(\|\tilde{e}_{s}\|_{\mathcal{T}}))^{2}\label{eq:dot_U-1}
\end{align}
where
\begin{equation}
\sigma_{\min}=\min\{\sigma^1,\cdots,\sigma^N\}\label{eq:sigma-min}
\end{equation}
and we used the fact $\hat{\theta}\t\varsigma(q,\dot{q},\varpi,\mu)\leq \mu\|\hat{\theta}\|\gamma_{\varsigma}(\|\tilde{e}_{s}\|_{\mathcal{T}}) \leq\sigma_{\min}\mu\|\hat{\theta}\|^{2}/2+\mu(\gamma_{\varsigma}(\|\tilde{e}_{s}\|_{\mathcal{T}}))^{2}/(2\sigma_{\min})$.
Invoking comparison lemma for (\ref{eq:dot_U-1}) yields
\begin{align*}
U(\hat{\theta}) & \leq\kappa^{-\sigma_{\min}}(\mu)U(\hat{\theta}(t_{0}))\\
 & \quad+\frac{1}{2\sigma_{\min}^{2}}(\gamma_{\varsigma}(\|\tilde{e}_{s}\|_{\mathcal{T}}))^{2}\left(1-\kappa^{-\sigma_{\min}}(\mu)\right)
\end{align*}
where $\kappa^{-\sigma_{\min}}(\mu)$ is defined as  in \eqref{eq:kappa} with $\alpha(s)=s$.
Therefore, the upper bound of $\hat{\theta}$ is
\begin{equation}
\|\hat{\theta}\|\leq\|\hat{\theta}(t_{0})\|+\sigma_{\min}^{-1}\gamma_{\varsigma}(\|\tilde{e}_{s}\|_{\mathcal{T}}).\label{eq:bound_hat_theta}
\end{equation}
 By (\ref{eq:tilde_e_s-1}), (\ref{eq:bound_Omega}) and (\ref{eq:bound_hat_theta}),
$\pi(q,\dot{q},\varpi,\hat{\theta},\mu)$ satisfies
\begin{equation}
\|\pi(q,\dot{q},\varpi,\hat{\theta},\mu)\|\leq\mu^{-\iota+2}\gamma_{\pi}(\|\tilde{e}_{s}\|_{\mathcal{T}})\label{eq:bound_pi}
\end{equation}
where $\gamma_{\pi}(s)=\|\bar K_2 \|s+\gamma_{\Omega}(s)(\|\hat{\theta}(t_{0})\|+\sigma_{\min}^{-1}\gamma_{\varsigma}(s))$ with $\bar K_2 = \diag\{k_2^1,\cdots,k_2^N\}\otimes I_n$.

Finally, we prove the boundedness of $\hat{\theta}^{i}$-dynamics
for $i\in\mathcal{V}$. Since (\ref{eq:dot_hat_theta^i}) is a non-homogeneous
linear differential equation, we can have an analytic solution of
$\hat{\theta}(t)$.
Substituting it into (\ref{eq:dot_hat_theta^i})
yields
\begin{align}
\dot{\hat{\theta}}(t) & =-\mu\bar\sigma\kappa^{-\bar\sigma}(\mu)\hat{\theta}(t_{0})+\varsigma(t)\nonumber \\
 & \quad-\mu\bar\sigma\kappa^{-\bar \sigma}(\mu)\sint_{t_{0}}^{t}\kappa^{\bar\sigma}(\mu(\tau)) \varsigma(\tau)\mathrm{d}\tau\label{eq:dot_hat_theta-1}
\end{align}
where
$\varsigma(q(t),\dot{q}(t),\varpi(t),\mu)$ is simplified as
$\varsigma(t)$, and
$\kappa^{-\bar\sigma}(\mu)=\exp\left(-\sint_{t_0}^t \mu(\tau)\bar\sigma\mathrm d\tau\right)$ and  $\kappa^{\bar\sigma}(\mu(\tau))=\exp\left(\sint_{t_0}^\tau \mu(s)\bar\sigma\mathrm ds\right)$ are diagonal matrix functions. Since $\sigma^i\geq\tilde{b}$ for $i\in \mathcal V$, one has $\sigma_{\min}\geq \tilde b$. Let $\sigma_{\max}=\max\{\sigma^1,\cdots,\sigma^N\}$, then
\begin{align}
 &\|\mu\bar\sigma\kappa^{-\bar \sigma}(\mu)\hat{\theta}(t_{0})\|\notag\\
 &\leq \sigma_{\max} \|\hat\theta(t_0)\| \mu \kappa^{-\sigma_{\min}}(\mu)\notag \\
 &\leq \sigma_{\max} \|\hat\theta(t_0)\| b \exp\left( \sint_{t_0}^t \dot \mu(\tau) /\mu(\tau)\mathrm d\tau \right)\kappa^{-\sigma_{\min}}(\mu)\notag \\
 &\leq \sigma_{\max}\mu(t_{0})\|\hat{\theta}(t_{0})\|\kappa^{(\tilde{b}-\sigma_{\min})}(\mu)\notag\\
 &  \leq\sigma_{\max}\mu(t_{0})\|\hat{\theta}(t_{0})\|\label{eq:first-term}
\end{align}
holds for $t\in\mathcal{T}_{p}$.
According to the definition of $\mu(t)$ in Definition \ref{Def:K_T}, we have $(\mu(t))^{-1}:\mathcal T_p\mapsto (0,b^{-1}]$ being strictly decreasing to zero. By the definition of a left limit, for any given $\varepsilon>0$ (no matter how small), there always exists $\delta>0$ such that for $0<T+t_0 -t<\delta$, it holds that $(\mu(t))^{-1}<\varepsilon$. Thus for any pair  $(\varepsilon,\delta)$, it follows that $\mu(t)<\infty$ for $t\in[t_0,T+t_0-\delta]$ and $\mu(t)=\infty$ for $t\in(T+t_0-\delta,T+t_0)$.
Referring to (\ref{eq:bound_varsigma}), we observe that $\varsigma(t)$
and consequently  $\mu\bar\sigma\kappa^{-\bar \sigma}(\mu)\sint_{t_{0}}^{t}\kappa^{\bar\sigma}(\mu(\tau))\varsigma(\tau)\mathrm{d}\tau$
are bounded for $t\in [t_0,T+t_{0}-\delta)$. Therefore, combining this with (\ref{eq:first-term}), we conclude that
$\dot{\hat{\theta}}$ is bounded for $t\in [t_0,T+t_{0}-\delta)$. As for
$t\in(T+t_0-\delta,T+t_0)$, by arithmetic operations of ultimate limit and \textit{L'Hospital}'s
rule,
\begin{align*}
 & \lim_{t\to T+t_{0}}\left(\varsigma(t)-\mu\bar\sigma\kappa^{-\bar \sigma}(\mu)\sint_{t_{0}}^{t}\kappa^{\bar\sigma}(\mu(\tau))\varsigma(\tau)\mathrm{d}\tau\right)\\
  & =\lim_{t\to T+t_{0}}\varsigma(t)-\lim_{t\to T+t_{0}}\mu(t)\bar \sigma\\
 & \quad\times\lim_{t\to T+t_{0}}\kappa^{-\bar\sigma}(\mu(t))\sint_{t_{0}}^{t}\kappa^{\bar\sigma}(\mu(\tau))\varsigma(\tau)\mathrm{d}\tau\\
 & =\lim_{t\to T+t_{0}}\varsigma(t)-\lim_{t\to T+t_{0}}\mu\bar\sigma\kappa^{\bar\sigma}(\mu)\varsigma(t)/(\mu\kappa^{\bar\sigma}(\mu)\bar\sigma)\\
 & =\lim_{t\to T+t_{0}}\varsigma(t)-\lim_{t\to T+t_{0}}\varsigma(t)=0.
\end{align*}
Thus, $\dot{\hat{\theta}}(t)$ is bounded for $t\in\mathcal{T}_{p}$.
This completes the proof.

\section{Proof of Theorem \ref{lem:SG-Criterion}}
 According to (\ref{eq:PTISS-1}) and (\ref{eq:PTISS-2}),
one has
\begin{gather*}
\begin{gathered}\underline{\alpha}_{1}(\|\chi_{1}\|) \leq V_{1}(\chi_{1})\leq \overline{\alpha}_{1}(\|\chi_{1}\|),\\
\dot{V}_{1}\leq-\alpha_{1}(\mu)V_{1}+l_{1}\alpha_{2}(\mu)V_{2}(\chi_2)+\gamma_{1}(\mu)\epsilon_{1}(\|d_{1}\|),
\end{gathered}
\\
\begin{gathered}\underline{\alpha}_{2}(\|\chi_{2}\|)\leq V_{2}(\chi_{2})\leq \overline{\alpha}_{2}(\|\chi_{2}\|),\\
\dot{V}_{2}\leq-\alpha_{2}(\mu)V_{2}+l_{2}\alpha_{1}(\mu)V_{1}(\chi_1)+\gamma_{2}(\mu)\epsilon_{2}(\|d_{2}\|).
\end{gathered}
\end{gather*}
Since $l_{1}l_{2}<1$, by some simple algebra operations, one has
\[
l_{1}l_{2}<\sqrt{l_{1}l_{2}}<1.
\]
 Define the Lyapunov function candidate for $\chi$-dynamics as
\begin{equation}
V(\chi)=l_{2}V_{1}(\chi_{1})+\sqrt{l_{1}l_{2}}V_{2}(\chi_{2}).\label{eq:V_chi}
\end{equation}
 Then the time derivative along the trajectories of (\ref{eq:PTISS-1}) and
(\ref{eq:PTISS-2}) is
\begin{align}
\dot{V} & (\chi)\leq l_{2}(-\alpha_{1}(\mu)V_{1}+l_{1}\alpha_{2}(\mu)V_{2}+\gamma_{1}(\mu)\epsilon_{1}(\|d_{1}\|))\nonumber \\
 & \quad+\sqrt{l_{1}l_{2}}(-\alpha_{2}(\mu)V_{2}+l_{2}\alpha_{1}(\mu)V_{1}+\gamma_{2}(\mu)\epsilon_{2}(\|d_{2}\|))\nonumber \\
 & =-\left(1-\sqrt{l_{1}l_{2}}\right)\left(l_{2}\alpha_{1}(\mu)V_{1}(\chi_{1})+\sqrt{l_{1}l_{2}}\alpha_{2}(\mu)V_{2}(\chi_{2})\right)\nonumber \\
 & \quad+l_{2}\gamma_{1}(\mu)\epsilon_{1}(\|d_{1}\|)+\sqrt{l_{1}l_{2}}\gamma_{2}(\mu)\epsilon_{2}(\|d_{2}\|)\nonumber \\
 & \leq-\left(1-\sqrt{l_{1}l_{2}}\right)\min\{\alpha_{1}(\mu),\alpha_{2}(\mu)\}V(\chi)\nonumber \\
 & \quad+\max\{\gamma_{1}(\mu),\gamma_{2}(\mu)\}\left(l_{2}\epsilon_{1}(\|d_{1}\|)+\sqrt{l_{1}l_{2}}\epsilon_{2}(\|d_{2}\|)\right)\nonumber \\
 & \leq-\alpha(\mu)V(\chi)+\gamma(\mu)\epsilon(\|d\|)\label{eq:dot_V}
\end{align}
with
\[
\begin{gathered}\alpha(\mu)=\left(1-\sqrt{l_{1}l_{2}}\right)\min\{\alpha_{1}(\mu),\alpha_{2}(\mu)\},\\
\gamma(\mu)=\gamma_{1}(\mu)+\gamma_{2}(\mu),\\
\epsilon(\|d\|)=l_{2}\epsilon_{1}(\|d\|)+\sqrt{l_{1}l_{2}}\epsilon_{2}(\|d\|)
\end{gathered}
\]
 where $\alpha$, $\gamma$, $\epsilon\in\mathcal{K}_{\infty}$, $d=[d_1\t,d_2\t]\t$,
 and we used the fact $\|d_i\|\leq \|d\|$ for $i=1,2$. We note that $\alpha(\mu)$
may be piecewise continuous but is integrable. Since $d_{i}$ belongs
to a compact set $\mathbb{D}_{i}$, by (\ref{eq:O_s}),
there exists a $\mathcal{K}_{\infty}$ function $\tilde{\epsilon}$
such that
\begin{align*}
\frac{\gamma(\mu)\epsilon(\|d\|)}{\alpha(\mu)}\leq & \left(\sup_{\mu\in\mathbb{R}_{p}}\frac{\gamma(\mu)}{\alpha(\mu)}\right)\epsilon(\|d\|_{\mathcal{T}})=\tilde{\epsilon}(\|d\|_{\mathcal{T}})<\infty,
\end{align*}
where $\|d\|_{\mathcal{T}}=\sup_{t\in\mathcal{T}_{p}}\|d(t)\|$. Invoking
comparison lemma for (\ref{eq:dot_V}) yields
 \begin{align}
V(\chi)
&\leq\kappa^{-1}(\alpha(\mu))V(\chi_{0})+ \int_{t_{0}}^{t}\exp\left(-\sint_{\tau}^{t}\alpha(\mu(s))\mathrm{d}s\right)\nonumber \\
 & \quad\times\gamma(\mu(\tau))\epsilon(\|d(\tau)\|)\mathrm{d}\tau\nonumber \\
 & = \kappa^{-1}(\alpha(\mu))V(\chi_{0})+ \kappa^{-1}(\alpha(\mu))\int_{t_0}^t\kappa(\alpha(\mu(\tau)))\notag\\
 &\quad  \times\gamma(\mu(\tau))\epsilon(\|d(\tau)\|)\mathrm{d}\tau \notag\\
 &  \leq \kappa^{-1}(\alpha(\mu))V(\chi_{0})+ \tilde \epsilon (\|d\|_{\mathcal T})\kappa^{-1}(\alpha(\mu))\notag\\
 &\quad  \times \int_{t_0}^t \kappa(\alpha(\mu(\tau)))\alpha(\mu(\tau)) \mathrm d\tau\notag\\
 &= \kappa^{-1}(\alpha(\mu))V(\chi_{0})+ \tilde \epsilon (\|d\|_{\mathcal T})\kappa^{-1}(\alpha(\mu))\notag\\
 &\quad \times \int_{0}^{\int_{t_0}^t\alpha(\mu(s))\mathrm ds}\kappa(\alpha(\mu(\tau)))\mathrm d \left(\sint_{t_0}^\tau \alpha(\mu(s))\mathrm ds \right)\notag\\
 &= \kappa^{-1}(\alpha(\mu))V(\chi_{0})+ \tilde \epsilon (\|d\|_{\mathcal T})(1-\kappa^{-1}(\alpha(\mu))). \label{eq:bound_V_chi}
 \end{align}
By (\ref{eq:PTISS-1}),  (\ref{eq:PTISS-2}) and \eqref{eq:bound_V_chi}, for $i=1,2$, $\|\chi_i\|$ satisfies
\begin{equation*}
\|\chi_i \|\leq \underline {\alpha}_i^{-1}\circ \overline \alpha(\|\chi(t_0)\|)
\end{equation*}
where $\overline \alpha(s)\in \mathcal K_\infty^e$ denoted as
\begin{align}
\overline \alpha(s) & = \max\{ l_2^{-1}, (l_1l_2)^{-1/2}\}\big(\tilde\epsilon(\|d\|_{\mathcal T})\notag\\
& \quad + \max\{l_2, \sqrt{l_1l_2}\}\left(\overline\alpha_1(s)+\overline\alpha_2(s)\right)\big).
\end{align}
Upon using the fact that $\left(\ssum_{i=1}^n|x_i|\right)^p\leq \ssum_{i=1}^n|x_i|^p$ for $x_i\in \mathbb R, i=1,\cdots,n$, $0<p\leq 1$, we then have
\begin{align}
\|\chi\| &= \sqrt{\chi_1\t\chi_1+\chi_2\t\chi_2} \leq \|\chi_1\|+\|\chi_2\|\notag\\
&\leq \ssum_{i=1}^2 \underline {\alpha}_i^{-1} \circ \overline \alpha(\|\chi(t_0)\|)\label{eq:bound_chi-1}.
\end{align}
 Therefore, \eqref{eq:bound_chi-1}
leads to (\ref{eq:bound_chi}) with $\tilde{\alpha}$ denoted as
\[
\tilde{\alpha}(s)=\ssum_{i=1}^2\underline{\alpha}_i^{-1}\circ \overline \alpha(s).
\]
 This completes the proof.
 
 \section{Proof of Lemma \ref{the:1}}
Define $\bar{R}=[r_{1},r_{2}]\otimes I_{n}\in\mathbb{R}^{nN\times nN}$
an orthogonal matrix with $r_{1}\in\mathbb R^N$ and $r_{2}\in \mathbb R^{N\times N-1}$ satisfy $r_{1}=1_{N}/\sqrt{N}$,
$r_{1}\t r_{2}=0$, $r_{2}\t r_{2}=I_{N-1}$ and $r_{2}r_{2}\t=\Pi_{N}$,
where $\Pi_{N}=I_{N}-\frac{1}{N}1_{N}1_{N}\t$. We introduce the following
state transformations
\begin{equation}
\xi=\mu^{\iota}\bar{R}\t e_{\varpi},\quad\psi=\mu^{\iota}\bar{R}\t e_{v},\quad\phi=\mu^{\iota}\bar{R}\t e_{y}\label{eq:state-trans}
\end{equation}
 where $\xi=[\xi_{1}\t,\xi_{2}\t]\t$, $\psi=[\psi_{1}\t,\psi_{2}\t]\t$
and $\phi=[\phi_{1}\t,\phi_{2}\t]\t$ with $\xi_1,\psi_1,\phi_1 \in \mathbb R^n$ and $\xi_2,\psi_2,\phi_2\in\mathbb R^{(N-1)n}$.
 By (\ref{eq:dot_varpi_i})
and (\ref{eq:dot_v_i}), for $t\in\mathcal{T}_{p}$, we have the following
derivatives
\begin{eqnarray}
\dot{\xi}_{1} & =&-c\mu^{\iota+1}\bar{r}_{1}\t\big( \nabla\tilde{F}(\varpi,\varpi^*)+ \nabla\tilde{F}(y,\varpi)\big)+\iota\tilde{\mu}\xi_{1},\label{eq:dot_xi_1}\\
\dot{\xi}_{2} & =&-c\mu^{\iota+1}\bar{r}_{2}\t\big( \nabla\tilde{F}(\varpi,\varpi^*)
+ \nabla\tilde{F}(y,\varpi)\big)-c\mu\psi_{2}\nonumber \\
 & & \quad-c\mu\bar{\mathcal{L}}_{r}\left(\xi_{2}+\phi_{2}\right)+\iota\tilde{\mu}\xi_{2},\label{eq:dot_xi_2}\\
\dot{\psi}_{1} & =&0,\quad\dot{\psi}_{2}=c\mu\bar{\mathcal{L}}_{r}(\xi_{2}+\phi_{2})+\iota\tilde{\mu}\psi_{2}\label{eq:dot_psi_1}
\end{eqnarray}
 where $\bar{r}_{1}=r_{1}\otimes I_{n}$, $\bar{r}_{2}=r_{2}\otimes I_{n}$,
$\bar{\mathcal{L}}_{r}=r_{2}\t\mathcal{L}r_{2}\otimes I_{n}$ and
$\tilde{\mu}=\dot{\mu}/\mu$. Define the Lyapunov function candidate
for $\xi$- and $\psi$-dynamics as $V(\xi,\psi)=V_{1}(\xi,\psi)+V_{2}(\xi,\psi)$
with
\begin{gather}
\begin{gathered}V_{1}(\xi,\psi)=\frac{\delta}{2}(\xi\t\xi+\psi\t\tilde{\mathcal{L}}_{r}^{-1}\psi),\\
V_{2}(\xi,\psi)=(\xi+\psi)\t(\xi+\psi)
\end{gathered}
\label{eq:V}
\end{gather}
 where $\delta$ is denoted in \eqref{eq:delta}, and
 $\tilde{\mathcal{L}}_{r}=\mbox{diag}\{I_{n},\bar{\mathcal{L}}_{r}\}\in\mathbb{R}^{nN\times nN}$
is a positive definite matrix. By (\ref{eq:dot_xi_1}),
(\ref{eq:dot_xi_2}) and (\ref{eq:dot_psi_1}), the time derivative of $V_1 (\xi,\psi)$ is
\begin{align}
&\dot{V}_{1} (\xi,\psi)\notag\\ &=\delta(\xi_{1}\t\dot{\xi}_{1}+\xi_{2}\t\dot{\xi}_{2}+\psi_{2}\t\bar{\mathcal{L}}_{r}^{-1}\dot{\psi}_{2})\nonumber \\
 & =\delta\iota\tilde{\mu}\left(\|\xi\|^{2}+\psi_{2}\t\bar{\mathcal{L}}_{r}^{-1}\psi_{2}\right)+\delta c\mu^{\iota+1}\big(-\xi_{1}\t\bar{r}_{1}\t
  \nabla\tilde{F}(\varpi,\varpi^*)\nonumber \\
 & \quad-\xi_{1}\t\bar{r}_{1}\t
  \nabla\tilde{F}(y,\varpi)-\xi_{2}\t\bar{r}_{2}\t \nabla\tilde{F}(\varpi,\varpi^*)-\xi_{2}\t\bar{r}_{2}\t \nabla\tilde{F}(y,\varpi)\big)\nonumber \\
 & \quad+\delta c\mu\left(-\xi_{2}\t\bar{\mathcal{L}}_{r}\xi_{2}-\xi_{2}\t\bar{\mathcal{L}}_{r}\phi_{2} +\psi_{2}\t\phi_{2}\right)\label{eq:dot_V_1}
\end{align}
where we used $\psi_{1}\t\dot{\psi}_{1}=0$.
Upon using (\ref{eq:bound_mu_derivative}), Young's inequality,
and Assumption \ref{ass:cost_func}, a few facts are $\tilde{\mu}\big(\|\xi\|^{2}+\psi_{2}\t\bar{\mathcal{L}}_{r}^{-1}\psi_{2}\big)\leq\tilde{b}\mu\big(\|\xi\|^{2}+\lambda_{2}^{-1}\|\psi\|^{2}\big)$,
$-\mu^{\iota}\big(\xi_{1}\t\bar{r}_{1}\t
 \nabla\tilde{F}(\varpi,\varpi^*)
+\xi_{2}\t\bar{r}_{2}\t \nabla\tilde{F}(\varpi,\varpi^*)\big)\leq-\rho_{c}\|\xi\|^{2}$,
$-\mu^{\iota}\big(\xi_{1}\t\bar{r}_{1}\t \nabla\tilde{F}(y,\varpi)+\xi_{2}\t\bar{r}_{2}\t \nabla\tilde{F}(y,\varpi)\big)\leq\rho_{c}\|\xi\|^{2}/2+\varrho_{c}^{2}\|\phi\|^{2}/(2\rho_{c})$,
$-\xi_{2}\t\bar{\mathcal{L}}_{r}\phi_{2}\leq\lambda_{2}\|\xi_{2}\|^{2}/2+\lambda_{N}^{2}\|\phi_{2}\|^{2}/(2\lambda_{2})$
and $\psi_{2}\t\phi_{2}\leq\|\psi_{2}\|^{2}/(8\delta)+2\delta\|\phi_{2}\|^{2}$.
Therefore, according to the above inequalities, $\dot{V}_{1}(\xi,\psi)$ satisfies
\begin{align}
\dot{V}_{1}(\xi,\psi) & \leq\big(-\rho_{c}c/2+\iota\tilde{b}\big)\delta\mu\|\xi\|_{2}-\lambda_{2}c\delta\mu\|\xi_{2}\|^{2}/2\nonumber \\
&\quad+\left(\varrho_{c}^{2}/(2\rho_{c})+\lambda_{N}^{2}/(2\lambda_{2}) +2\delta\right)c\delta\mu\|\phi\|^{2}\notag\\
 & \quad+\big(\iota\tilde{b}/\lambda_{2}+c/(8\delta)\big)\delta\mu\|\psi\|^{2}\label{eq:dot_V1}
\end{align}
 where we used $\|\phi_{2}\|\leq\|\phi\|$ and $\|\psi_{2}\|\leq\|\psi\|$.
The derivative of $V_{2}(\xi,\psi)$ can be expressed as
\begin{align*}
&\dot{V}_{2}(\xi,\psi)\notag\\ &=\xi_{1}\t\dot{\xi}_{1}+\xi_{2}\t\dot{\xi}_{2}+\xi_{2}\t\dot{\psi}_{2}+\psi_{1}\t\dot{\xi}_{1}+\psi_{2}\t\dot{\xi}_{2}+\psi_{2}\t\dot{\psi}_{2}\\
&=  -c\mu^{\iota+1}\left(\xi_{1}\t\bar{r}_{1}\t+\xi_{2}\t\bar{r}_{2}\t\right)\big( \nabla\tilde{F}(\varpi,\varpi^*)+\nabla\tilde{F}(y,\varpi)\big)\\
 & \quad-c\mu^{\iota+1}\left(\psi_{1}\t\bar{r}_{1}\t+\psi_{2}\t\bar{r}_{2}\t\right)\big( \nabla\tilde{F}(\varpi,\varpi^*)+ \nabla\tilde{F}(y,\varpi)\big)\\
 & \quad -c\mu\left(\xi_{2}\t\psi_{2}+\psi_{2}\t\psi_{2}\right)+\iota\tilde{\mu}\|\xi+\psi\|^{2}.
\end{align*}
 By Young's inequality and (\ref{eq:bound_mu_derivative}), one has
$-\mu^{\iota}(\psi_{1}\t\bar{r}_{1}\t+\psi_{2}\t\bar{r}_{2}\t) \nabla\tilde{F}(\varpi,\varpi^*)\leq\|\psi\|^{2}/8+2\varrho_{c}^{2}\|\xi\|^{2}$,
$-\mu^{\iota}(\psi_{1}\t\bar{r}_{1}\t+\psi_{2}\t\bar{r}_{2}\t) \nabla\tilde{F}(y,\varpi)\leq\|\psi\|^{2}/8+2\varrho_{c}^{2}\|\phi\|^{2}$,
$-\xi_{2}\t\psi_{2}\leq2\|\xi_{2}\|^{2}+\|\psi_{2}\|^{2}/8$ and $\iota\tilde{\mu}\|\xi+\psi\|^{2}\leq2\iota\tilde{b}\mu\left(\|\xi\|^{2}+\|\psi\|^{2}\right)$.
Then by the inequalities under (\ref{eq:dot_V_1}), and the upper
inequalities, $\dot{V}_{2}(\xi,\psi)$ satisfies
\begin{align}
&\dot{V}_{2}(\xi,\psi) \notag\\
& \leq\big(-\rho_{c}c/2+2c\varrho_{c}^{2}+2\iota\tilde{b}\big)\mu\|\xi\|^{2} +\big(5c/8+2\iota\tilde{b}\big)\mu\|\psi\|^{2}\nonumber \\
 & \quad+2c\mu\|\xi_{2}\|^{2}+c\left(\varrho_{c}^{2}/(2\rho_{c})+2\varrho_{c}^{2}\right)\mu\|\phi\|^{2}.\label{eq:dot_V2}
\end{align}
 Therefore, by (\ref{eq:dot_V1}) and (\ref{eq:dot_V2}), $\dot{V}(\xi,\psi)$
satisfies
\begin{align*}
\dot{V}(\xi,\psi)\leq & \big(-\rho_{c}c\delta/2-\rho_{c}c/2+2c\varrho_{c}^{2}+\delta\iota\tilde{b}+2\iota\tilde{b}\big)\mu\|\xi\|^{2}\\
 & +\big(-c/2+\delta\iota\tilde{b}/\lambda_{2}+2\iota\tilde{b}\big)\mu\|\psi\|^{2}\\
 & +c_{\Delta}\mu\|\phi\|^{2}+\left(-\delta\lambda_{2}c/2+2c\right)\mu\|\xi_{2}\|^{2}
\end{align*}
where $c_\Delta$ is denoted in \eqref{eq:c-delta}.
Recall the definition of $\delta$ in (\ref{eq:delta}), we have $-\rho_{c}\delta c/2-\rho_{c}c/2+2c\varrho_{c}^{2}\leq-c/2$
and $-\delta\lambda_{2}c/2+2c\leq0$. Therefore, for $c$ denoted
in (\ref{eq:c}), $\dot{V}(\xi,\psi)$ further satisfies
\begin{equation}
\dot{V}(\xi,\psi)\leq-\bar{\delta}c^{*}\mu\left(\|\xi\|^{2}+\|\psi\|^{2}\right)+c_{\Delta}\mu\|\phi\|^{2}.\label{eq:dot_V-1}
\end{equation}
We define the Lyapunov functions for $\tilde{e}_{r}$- and $\tilde{e}_{s}'$-dynamics
as
\begin{gather}
U(\tilde{e}_{r})=\frac{\delta}{2}\tilde{e}_{r}\t\Lambda_{r,1}\tilde{e}_{r}+\frac{1}{2}\tilde{e}_{r}\t\Lambda_{r,2}\tilde{e}_{r},\label{eq:U}\\
W(\tilde{e}_{s}')=\left(\tilde{e}_{s}'\right)\t\Lambda_{s}\tilde{e}_{s}'\label{eq:W}
\end{gather}
where $\Lambda_{r,1}=\mbox{diag}\{I_{mn},\bar{R}\tilde{\mathcal{L}}_{r}^{-1}\bar{R}\t\}$, $\Lambda_{r,2}=[I_{mn}\;I_{mn}]\t [I_{mn}\; I_{mn}]$ and $\Lambda_{s}=\mbox{diag}\{I_{nN},M(q),I_{pN}\}$.
We note
\[
U(\tilde{e}_{r})=V(\xi,\psi),\quad\dot{U}(\tilde{e}_{r})=\dot{V}(\xi,\psi).
\]
We decompose $\tilde{e}_{s}$ in (\ref{eq:tilde_e_r})
as
\begin{equation}
\tilde{e}_{s}=\left[\tilde{e}_{s,1}\t,\tilde{e}_{s,2}\t\right]\t\notag
\end{equation}
with $\tilde{e}_{s,1}=\mu^{\iota}e_{y}$ and $\tilde{e}_{s,2}=\mu^{\iota-1}\left(\dot{q}+\mu\bar{K}_{1}e_{y}\right)$.
Since $\bar{R}$ in \eqref{eq:state-trans} is an orthogonal matrix, by (\ref{eq:tilde_e_r}),
one has $\|\phi\|\leq\|\tilde{e}_{s,1}\|$, then
according to (\ref{eq:dot_V-1}), (\ref{eq:U}) and the fact $\|\tilde{e}_{s,1}\|^{2}\leq W(\tilde{e}_{s}')$,
\begin{equation}
\begin{gathered}\underline{\delta}\|\tilde{e}_{r}\|^{2}\leq U(\tilde{e}_{r})\leq\bar{\delta}\|\tilde{e}_{r}\|^{2},\\
\dot{U}(\tilde{e}_{r})\leq-c^{*}\mu U(\tilde{e}_{r})+c_{\Delta}\mu W(\tilde{e}_{s}').
\end{gathered}
\label{eq:dot_U}
\end{equation}

Now we prove $\tilde{e}_{s}'$-dynamics admits a prescribed-time ISS
Lyapunov function.
It's worth noting that $\Lambda_{s}$ in \eqref{eq:W} is positive definition according to Property
\ref{prop:1-3}. Then  $2\tilde{e}_{s,1}\t\dot \tilde{e}_{s,1}$ can be expressed as
\begin{align*}
2\tilde{e}_{s,1}\t\dot{\tilde{e}}_{s,1} & =2\tilde{e}_{s,1}\t(\iota\mu^{\iota-1}\dot{\mu}e_{y}+\mu^{\iota}(\dot{q}-\dot{\varpi}))\\
 & =2\tilde{e}_{s,1}\t(\iota\tilde{\mu}\tilde{e}_{s,1}+\mu\tilde{e}_{s,2} -\mu\bar{K}_{1}\tilde{e}_{s,1}-\mu^{\iota}\dot{\varpi})
\end{align*}
where we used $\mu\dot{q}=\mu\tilde{e}_{s,2}-\mu^{\iota+1}\bar{K}_{1}e_{y}$
according to (\ref{eq:tilde_e_s-1}). By Young's inequality, (\ref{eq:bound_mu_derivative})
and (\ref{eq:k-1}), one has $2\tilde{\mu}\iota\|\tilde{e}_{s,1}\|^{2}\leq2\tilde{b}\iota\mu\|\tilde{e}_{s,1}\|^{2}$,
$2\mu\tilde{e}_{s,1}\t\tilde{e}_{s,2}\leq\mu\|\tilde{e}_{s,1}\|^{2}+\mu\|\tilde{e}_{s,1}\|^{2}$,
$-2\mu^{l}\tilde{e}_{s,1}\t\dot{\varpi}\leq\mu\|\tilde{e}_{s,1}\|^{2}+\mu^{2\iota-1}\|\dot{\varpi}\|^{2}$
and $-2\mu\tilde{e}_{s,1}\t\bar{K}_{1}\tilde{e}_{s,1}=-2\mu\tilde{e}_{s,1}\t\bar{K}_{1}^{*}\tilde{e}_{s,1}-2(\tilde{b}\iota+1)\mu\|\tilde{e}_{s,1}\|^{2}$,
where $\bar{K}_{1}^{*}=\diag\{k_{1}^{1,*},\cdots,k_{1}^{N,*}\}\otimes I_{n}$.
Substituting the above inequalities into $2\tilde{e}_{s}\t\dot{\tilde{e}}_{s}$
yields
\begin{eqnarray}
2\tilde{e}_{s,1}\t\dot{\tilde{e}}_{s,1}\leq-2\mu\tilde{e}_{s,1}\t\bar{K}_{1}^{*}\tilde{e}_{s,1}+\mu\|\tilde{e}_{s,2}\|^{2}+\mu^{2\iota-1}\|\dot{\varpi}\|^{2}.\label{eq:dot_W_1}
\end{eqnarray}
We note
\[
\frac{\mathrm{d}\tilde{e}\t M(q)\tilde{e}_{s,2}}{\mathrm{d}t}=2\tilde{e}_{s,2}\t M(q)\dot{\tilde{e}}_{s,2}+\tilde{e}_{s,2}\t\dot{M}(q)\tilde{e}_{s,2}.
\]
 According to (\ref{eq:Euler-Lag}), (\ref{eq:tau^i}) and (\ref{eq:tilde_e_s-1}),
one has
\begin{align}
 & \tilde{e}_{s,2}\t M(q)\dot{\tilde{e}}_{s,2}\nonumber \\
 & =\tilde{e}_{s,2}\t\Big[(\iota-1)\tilde{\mu}M(q)\tilde{e}_{s,2}+\tilde{\mu}M(q)\bar{K}_{1}\tilde{e}_{s,1}\nonumber \\
 & \quad+\mu^{\iota}M(q)\bar{K}_{1}(\dot{q}-\dot{\varpi})+\mu^{\iota-1}(\pi(q,\dot{q},\varpi,\hat{\theta},\mu)-C(q,\dot{q})\dot{q})\Big]\nonumber \\
 & =\tilde{e}_{s,2}\t\mu^{\iota-1}\Big[\pi(q,\dot{q},\varpi,\hat{\theta},\mu)+M(q)z_{1}+C(q,\dot{q})z_{2}\nonumber \\
 & \quad-\mu M(q)\bar{K}_{1}\dot{\varpi}\Big]-\tilde{e}_{s,2}\t C(q,\dot{q})\tilde{e}_{s,2}\label{eq:dot_tilde_e_s2}
\end{align}
 where $z_{1}$ and $z_{2}$ are defined in (\ref{eq:z_1^i}) and (\ref{eq:z_2^i}), and we used the fact $\dot{q}=\mu^{-\iota+1}\left(\tilde{e}_{s,2}-\bar{K}_{1}\tilde{e}_{s,1}\right)$.
Note that $\tilde{e}_{s,2}$ can be decomposed as $\tilde{e}_{s,2}=\left[(\tilde{e}_{s,2}^{1})\t,\cdots,(\tilde{e}_{s,2}^{N})\t\right]\t$
with $\tilde{e}_{s,2}^{i}=\mu^{\iota-1}(q^{i}+\mu k_{1}^{i}e_{y}^{i})$
for $i\in\mathcal{V}$. By Young's inequality, one has
\begin{align}
 & -2\mu^{\iota}\tilde{e}_{s,2}\t M(q)\bar{K}_{1}\dot{\varpi}\nonumber \\
 & =\ssum_{i=1}^{N}\left[(-2\mu^{\iota})k_{1}^{i}(\tilde{e}_{s,2}^{i})\t M^{i}(q^{i})\dot{\varpi}^{i}\right]\nonumber \\
 & \leq\ssum_{i=1}^{N}\left[\mu(k_{1}^{i})^{2}(k_{\overline{m}}^{i})^{2}\|\tilde{e}_{s,2}^{i}\|^{2}+\mu^{2\iota-1}\|\dot{\varpi}^{i}\|^{2}\right]\nonumber \\
 & =\mu\tilde{e}_{s,2}\t\bar{K}_{1}^{2} \bar{K}_{\overline{m}}^{2}\tilde{e}_{s,2}+\mu^{2\iota-1}\|\dot{\varpi}\|^{2}\label{eq:term_dot_varpi}
\end{align}
 where $\bar{K}_{\overline{m}}=\diag\{k_{\overline{m}}^{1},\cdots,k_{\overline{m}}^{N}\}\otimes I_{n}$
and we note that $\bar{K}_{1}^{2}\bar{K}_{\overline{m}}^{2}$ is a positive
definite matrix. Taking time derivative of $\tilde{\theta}\t\tilde{\theta}$
yields
\[
\tilde{\theta}\t\dot{\tilde{\theta}}=-\tilde{\theta}\t\varsigma(q,\dot{q},\varpi^{i},\mu) +\mu\tilde{\theta}\t\bar{\sigma}\hat{\theta}
\]
where $\varsigma(\cdot)$ is defined in (\ref{eq:varsigma^i}). According
to Property \ref{prop:1-1} and \ref{prop:1-2}, one has $-2\tilde{e}_{s,2}\t C(q,\dot{q})\tilde{e}_{s,2}+\tilde{e}_{s,2}\t\dot{M}(q)\tilde{e}_{s,2}=0$
and $M(q)z_{1}+C(q,\dot{q})z_{2}=\Omega(q,\dot{q},z_{1},z_{2})\theta$.
Substituting $\pi(q,\dot{q},\varpi,\hat{\theta},\mu)$ in (\ref{eq:tau^i})
and (\ref{eq:term_dot_varpi}) into (\ref{eq:dot_tilde_e_s2}), we
have
\begin{align}
 & 2\tilde{e}_{s,2}\t M(q)\dot{\tilde{e}}_{s,2}+\tilde{e}_{s,2}\t\dot{M}(q)\tilde{e}_{s,2}+2\tilde{\theta}\t\dot{\tilde{\theta}}\nonumber \\
 & =-\mu\tilde{e}_{s,2}\t\left(2\bar{K}_{2}-\bar{K}_{1}\bar{K}_{\bar{m}}\right)\tilde{e}_{s,2}+2\mu\tilde{\theta}\t\bar{\sigma}\hat{\theta}+\mu^{2\iota-1}\|\dot{\varpi}\|^{2}\nonumber \\
 & \leq-2\mu\tilde{e}_{s,2}\t\bar{K}_{2}^{*}\tilde{e}_{s,2}-\mu\|\tilde{e}_{s,2}\|^{2}-\mu\tilde{\theta}\t\bar{\sigma}\tilde{\theta}\nonumber \\
 & \quad+\mu\theta\t\bar{\sigma}\theta+\mu^{2\iota-1}\|\dot{\varpi}\|^{2}\label{eq:dot_W_2}
\end{align}
where $\bar{K}_{2}^{*}=\diag\{k_{2}^{1,*},\cdots,k_{2}^{N,*}\}\otimes I_{n}$ with $k_2^{i,*}$ been denoted in \eqref{eq:k-1} for $i\in \mathcal V$.
We note
\begin{gather*}
 \underline \varepsilon\|\tilde{e}_{s}'\|^{2}\leq W(\tilde{e}_{s}')\leq \overline \varepsilon \|\tilde{e}_{s}'\|^{2},\\
\dot{W}(\tilde{e}_{s}')=2\tilde{e}_{s,1}\t\dot{\tilde{e}}_{s,1}+2\tilde{e}_{s,2}\t M(q)\dot{\tilde{e}}_{s,2}+\tilde{e}_{s,2}\t\dot{M}(q)\tilde{e}_{s,2}+2\tilde{\theta}\t\dot{\tilde{\theta}}
\end{gather*}
where $\underline \varepsilon$ and $\overline \varepsilon$ are denoted in \eqref{eq:delta}.
 By (\ref{eq:bound_dot_e_r}) and the fact $\|\dot{\varpi}\|\leq\|\dot{e}_{r}\|$,
one has
\[
2\mu^{2\iota-1}\|\dot{\varpi}\|^{2}\leq c_{s}\mu(\|\tilde{e}_{r}\|^{2}+\|\tilde{e}_{s,1}\|^{2})
\]
where $c_s$ is denoted in \eqref{eq:c-delta}.
Then by (\ref{eq:W}),(\ref{eq:dot_W_1}), and (\ref{eq:dot_W_2}),
$\dot{W}(\tilde{e}_{s}')$ further satisfies
\begin{eqnarray}
\dot{W}(\tilde{e}_{s}')&\leq & -\tilde{k}_{1}\mu\|\tilde{e}_{s,1}\|^{2}-\tilde{k}_{2}\mu\tilde{e}_{s,2}\t M(q)\tilde{e}_{s,2}\nonumber \\
 & & -\sigma_{\min}\mu\|\tilde{\theta}\|^{2}+\sigma_{\max}\mu\|\theta\|^{2}+c_{s}\mu\|\tilde{e}_{r}\|^{2}\nonumber \\
&\leq & -\tilde{k}\mu W(\tilde{e}_{s}')+\sigma_{\max}\mu\|\theta\|^{2}+(c_{s}/\underline{\delta})\mu U(\tilde{e}_{r})\label{eq:dot_W}
\end{eqnarray}
where $\tilde k_1$, $\tilde k_2$ and $\tilde k$ have been denoted under \eqref{eq:con2}, and we used the fact
\begin{align*}
-2\mu\tilde{e}_{s,2}\t\bar{K}_{2}^{*}\tilde{e}_{s,2} & =-2\mu\ssum_{i=1}^{N}k_{2}^{i}(\tilde{e}_{s,2}^{i})\t\tilde{e}_{s,2}^{i}\\
 & \leq-2\mu\ssum_{i=1}^{N}k_{2}^{i}(k_{\bar{m}}^{i})^{-1}(\tilde{e}_{s,2}^{i})\t M^{i}(q^{i})\tilde{e}_{s,2}^{i}\\
 & =2\mu\tilde{k}_{2}\tilde{e}_{s,2}\t M(q)\tilde{e}_{s,2}.
\end{align*}
Therefore, (\ref{eq:dot_U}) and (\ref{eq:dot_W}) lead to (\ref{eq:con1})
and (\ref{eq:con2}). This completes the proof.

 }
\end{document}